\documentclass[12pt]{amsart}
\usepackage{amsthm}
\usepackage{amsmath, mathrsfs}
\usepackage{amssymb}
\usepackage{hyperref}
\usepackage{enumerate}
\usepackage{latexsym,array}
\usepackage{amsfonts}
\usepackage{shadow}

\newcommand{\bigR}{\mathcal{R}}
\newcommand{\bigS}{\mathcal{S}}

\newcommand{\w}{\omega}
\newtheorem{Pa}{Paper}[section]
\newtheorem{Tm}[Pa]{{\bf Theorem}}

\newtheorem{Dn}[Pa]{{\bf Definition}}
\newtheorem{Cy}[Pa]{{\bf Corollary}}

\newtheorem{Rk}[Pa]{{\bf Remark}}
\newtheorem{Pn}[Pa]{{\bf Proposition}}
\newtheorem{Pb}[Pa]{{\bf Problem}}

\date{}
\author{Daniel Alpay}
\address{(DA) Department of Mathematics
\newline
Ben Gurion University of the Negev \newline P.O.B. 653,
\newline
Be'er Sheva 84105, \newline ISRAEL}
\email{dany@math.bgu.ac.il}

\thanks{D. Alpay thanks the
Earl Katz family for endowing the chair which supported his
research, and the Binational Science Foundation Grant number
2010117. 
It is a pleasure to thank the referee for her/his
careful reading of the manuscript.
}
\author{Guy Salomon}
\address{(GS) Department of Mathematics
\newline
Ben Gurion University of the Negev \newline P.O.B. 653,
\newline
Be'er Sheva 84105, \newline ISRAEL} \email{guysal@math.bgu.ac.il}
\title[New Topological $\mathbb C$-Algebras]
{New Topological $\mathbb C$-Algebras with Applications in Linear
Systems Theory}
\begin{document}
\maketitle

\tableofcontents
\begin{abstract}
Motivated by the Schwartz space of tempered distributions
$\mathscr S^\prime$ and the Kondratiev space of stochastic
distributions $\mathcal S_{-1}$ we define a wide family of
nuclear spaces which are increasing unions of (duals of) Hilbert
spaces $\mathscr H_p^\prime,p\in\mathbb N$, with decreasing norms
$\|\cdot\|_{p}$. The elements of these spaces are functions on a
free commutative monoid. We characterize those rings in this
family which satisfy an inequality of the form $\|f *
g\|_{p} \leq A(p-q) \|f\|_{q}\|g\|_{p}$ for all $p\ge q+d$, where
$*$ denotes the convolution in the monoid, $A(p-q)$ is a
strictly positive number and $d$ is a fixed natural number (in
this case we obtain commutative topological $\mathbb C$-algebras). Such an
inequality holds in $\mathcal S_{-1}$, but not in $\mathscr
S^\prime$. We give an example of such a ring which contains
$\mathscr S^\prime$. We characterize invertible elements in these
rings and present applications to linear system theory.
\end{abstract}
\section{Introduction}
\setcounter{equation}{0} As is well known, the Schwartz space
$\mathscr S^\prime$ of complex tempered distributions   can be
viewed as the space of sequences of complex numbers $f=(f_n)_{n
\in \mathbb{N}_0}$ subject to
\[
\sum_{n \in \mathbb{N}_0} {
 (n+1)^{-2p}|f_n| ^2<\infty} \text{ for
some $p \in \mathbb{N}$}.
\]
Setting
\[
\|f\|_p=\left(\sum_{n \in \mathbb{N}_0}{(n+1)^{-2p}}|f_n|^2\right)^{1/2}<\infty,
\]
one can represent $\mathscr S^\prime$ as a union of an increasing
sequence of  Hilbert spaces $\mathscr H_1^\prime,\mathscr
H_2^\prime,\ldots $of complex sequences, with decreasing norms:
\begin{equation}
\label{eq:Hp} \mathscr H_p^\prime=\left\{f=(f_n)_{n \in
\mathbb{N}_0}:\, \|f\|_p<\infty\right\}.
\end{equation}

In Hida's white noise space theory, a counterpart of $\mathscr
S^\prime$ was introduced by Kondratiev, see \cite{MR1408433} and
the references therein, in the following way. We begin with a
definition:
\begin{Dn}
$\ell$ denotes the set of sequences of elements of $\mathbb N_0$,
indexed by $\mathbb N$,
\[
(\alpha_1,\alpha_2,\ldots)
\]
where $\alpha_j\not=0$ for at most a finite number of indices.
\label{def:ell}
\end{Dn}
The stochastic counterpart of $\mathscr S^\prime$ is the space
$\bigS_{-1}$ of families of complex numbers
$f=(f_\alpha)_{\alpha\in\ell}$ indexed by $\ell$ and such that
\[
\sum_{\alpha \in \ell} {|f_\alpha| ^2 (2 \mathbb{N})^{-\alpha
p}<\infty} \text{ for some  $p \in \mathbb{N}$}.
\]
In the above expression, one sets
\[
(2\mathbb{N})^{\alpha}=2^{\alpha_1}4^{\alpha_2}6^{\alpha_3}\cdots
\]
We here set
\begin{equation}
\label{eq:norms}
\|f\|_p=\|(f_\alpha)_{\alpha\in\ell}\|_p=\left(\sum_{\alpha \in \ell}
{|f_\alpha| ^2 (2 \mathbb{N})^{-\alpha p}}\right)^{1/2},
\end{equation}
and denote
\[
\mathscr K_p^\prime=\left\{f=(f_\alpha)_{\alpha \in \ell}:\,
\|f\|_p<\infty\right\}.
\]
In a way similar to $\mathscr S^\prime$, the space $\bigS_{-1}$
is the union of the increasing sequence of Hilbert spaces $\mathscr K_1^\prime,\mathscr K_2^\prime,\ldots$ with
decreasing norms \eqref{eq:norms}. The elements of $\bigS_{-1}$ are called
stochastic distributions and play an important role in stochastic
partial differential equations, see \cite{MR1408433}. We also
refer to the papers \cite{aa_goh,al_acap,alp} where $\bigS_{-1}$
is used to develop a new approach to linear stochastic systems.
Recall that the convolution of two elements of $\bigS_{-1}$,
$f=(f_\alpha)_{\alpha \in \ell}$ and $g=(g_\alpha)_{\alpha \in
\ell}$, (which is called the Wick product and denoted by
$\diamond$) is defined by
\[
f\diamond g=\left(\sum_{\beta \leq \alpha}f_\beta
g_{\alpha-\beta}\right)_{\alpha \in \ell}
\]
and satisfies the inequality
\begin{equation}
\label{eq:vage}
\|f \diamond g\|_{p} \leq A(p-q)
\|f\|_{q}\|g\|_{p},\quad \text{ for all } p\ge q+2,
\end{equation}
where
\[
A(p-q)=\left(\sum _{\alpha \in \ell} (2 \mathbb{N})^{-(p-q)
\alpha}\right)^{\frac12}
\]
is finite. This inequality is due to V\aa ge, see \cite{vage96},
\cite{MR1408433}. It expresses in particular that the
multiplication operator
\[
g\mapsto f\diamond g
\]
is bounded from the Hilbert space $\mathscr K_p^\prime$ into
itself where $f\in\mathscr K_q^\prime$ and $p\ge q+2$. It
plays a key role in the applications mentioned above. \\

In view of \eqref{eq:vage} it is a natural question to ask if a
similar inequality holds in $\mathscr S^\prime$. Here lies an
important structural difference between $\mathscr S^\prime$ and
$\bigS_{-1}$. If $f=(f_n)_{n\in\mathbb N_0}$ and
$g=(g_n)_{n\in\mathbb N_0}$ belong to $\mathscr S^\prime$, their
convolution which is defined by
\[
f * g=\left(\sum_{m \leq n}f_m g_{n-m}\right)_{n \in \mathbb N_0}
\]
also belongs to $\mathscr S^\prime$.
Nevertheless, as will be proved in the sequel (see Corollary
\ref{pn:Schwartz}), one cannot have an inequality of the kind
\eqref{eq:vage}, that is:
\begin{equation}
\label{vage123345}
\|f * g\|_{p} \leq B(p-q)\|f\|_{q}\|g\|_{p},\quad \text{ for all } p\ge q+d,
\end{equation}
where $d\in\mathbb N$ is preassigned, $B(p-q)>0$ is a constant
which depends only on $p-q$ in $\mathbb N$, and where $f$ runs
through $\mathscr H_q^\prime$ and $g$ runs through
$\mathscr H_p^\prime$.
Since such an inequality does not
hold, the origin of the present study was to find nuclear spaces
containing $\mathscr S^\prime$ such that an appropriate
inequality of the type \eqref{eq:vage} holds for the convolution.\\

More generally, in the present paper we define a wide family of
nuclear spaces in terms of positive functions over a commutative monoid, and give
a characterization of those in which an inequality of the type
\eqref{eq:vage} holds. Since such an inequality was first proved
by V\aa ge (in the setting of the Kondratiev space of stochastic
distributions) we call these spaces V{\aa}ge spaces. We show that
these spaces are in particular topological $\mathbb C$-algebras, and give a
characterization of their invertible elements. We then consider
the tensor product of two V{\aa}ge spaces,
and show that it is a V{\aa}ge space too.\\

The Schwartz space of tempered distributions $\mathscr S^\prime$
is not a V{\aa}ge space. We define a V{\aa}ge space,
containing $\mathscr S^\prime$. This space
is the dual of a space which is included in the
Schwartz space of test functions $\mathscr S$, consists of entire
functions, and is invariant under the Fourier transform.
One can thus define the Fourier transform on its dual, and study,
as suggested to us by Palle Jorgensen, connections with the theory of
hyperfunctions (see \cite{MR1311923} for the latter). This will be done in
another publication. We present
some important properties of this space, and characterize it both in
terms of sequences and in terms of entire functions. \\

The paper consists of seven sections besides the introduction, and
we now describe its content. A family of spaces of functions over a commutative monoid
which includes the space $\mathscr S^\prime$, and which we call
{\sl regular admissible spaces}, is introduced in Section
\ref{sec:2}. In Section \ref{sec:3} we characterize regular
admissible spaces in which convolution satisfies an inequality of
the type \eqref{eq:vage}. We also prove that these spaces are
topological $\mathbb C$-algebras, which we call V\aa ge spaces. Invertible
elements in these rings are characterized in Section \ref{sec:4}.
In Section \ref{sec:5} we prove that the tensor product of two
V\aa ge spaces is a V\aa ge space. The last three sections are
devoted to examples and applications. In Section
\ref{sec:schwartz} we define a V\aa ge space which contains
$\mathscr S^\prime$. Some results of E. Hille on Hermite series
play an important role in the arguments. In Section \ref{sec:7} we
consider the Kondratiev space. Finally, applications to linear
system theory are outlined in Section \ref{sec:8}.


\section{A new family of nuclear spaces and Gelfand triples}
\label{sec:2}
\setcounter{equation}{0}
In this section we introduce a family of nuclear spaces of
functions over a commutative monoid which we use in the sequel. We begin with a definition
and a preliminary result on such positive functions. Let $A$ be a
subset of $\mathbb{N}$. We denote
\begin{equation}
\label{eq:lA} \ell_A= \mathbb N_0^{(A)}=\left\{\alpha \in \mathbb N_0^A : \text{supp}(\alpha) \text{ is finite} \right\}=\oplus_{n\in A}\mathbb N_0 e_n,
\end{equation}
where, for $n\in\mathbb N$, we have denoted by $e_n$ the sequence
with all elements equal to 0, at exception of the $n$-th one,
equal to 1.
For two elements
\[\alpha=\sum_{n \in A} \alpha_n e_n\quad{\rm and}\quad
\beta=\sum_{n \in A} \beta_n e_n
\]
in $\ell_A$, we define $\alpha+\beta=\sum_{n \in A}
(\alpha_n+\beta_n) e_n$. In other words, $(\ell_A,+,0)$ is the
free commutative monoid generated by the countable (or finite)
set $\{e_n\}_{n\in A}$. Moreover, we consider the following
partial order induced by the addition above: For $\alpha,\beta
\in \ell_A$, we define $\alpha \leq \beta$ if there exists $\gamma
\in \ell_A$ such that $\alpha+\gamma=\beta$.\\

\begin{Dn}
\label{def:admissible}
Let $A$ be a subset of $\mathbb{N}$, and
let $\ell_A$ be defined by \eqref{eq:lA}. A positive function $a:\ell_a \to \mathbb R$ (that is $\alpha \mapsto a_\alpha$ where $a_\alpha>0$ for any $\alpha  \in \ell_A$) is called {\sl
admissible} if $a_0=1$ and $a_{e_n}>1$
for all $n \in A$.\\
Let $d\in\mathbb N$. The admissible function $a$ is called {\sl
$d$-regular} (or simply {regular}) if furthermore
\begin{equation}
\label{eq:ineq1} \sum _{n \in A} \frac 1{a_{e_n}^d-1}<
\infty.
\end{equation}
It is {\sl superexponential} (resp. {\sl exponential})  if
\begin{equation}
\label{eq:ineq}
a_\alpha a_\beta \leq a_{\alpha
+\beta}\quad \text{(resp. $a_\alpha a_\beta=a_{\alpha
+\beta}$)}\quad \forall\alpha,\beta\in\ell_A.
\end{equation}
\end{Dn}

Two examples of exponential regular
admissible positive functions are as follows:
\begin{enumerate}[(a)]
\item
The set $A$ has cardinal one. Then, $\ell_A=\mathbb N_0$.
We take $a_n=c^n$ with $c>1$.\\
\item
We set $A=\mathbb{N}$. Then, $\ell_A=\ell$, where $\ell$ is as in
Definition \ref{def:ell}. We take
\begin{equation}
\label{eq:zhang}
a_\alpha=(2
\mathbb{N})^{\alpha}=2^{\alpha_1}4^{\alpha_2}6^{\alpha_3}\cdots,\quad
\alpha\in\ell.
\end{equation}
\end{enumerate}
As mentioned in the introduction, this last example occurs in
Hida's white noise space theory, in the definition of spaces of
stochastic distributions.

\begin{Pn} \label{Pn:Converges}
Let $a:\ell_A \to \mathbb R$  be a superexponential $d$-regular
admissible positive function. Then,
\[\sum _{\alpha \in \ell_A}
a_\alpha^{-d} < \infty.
\]
Furthermore, if $a$ is exponential rather than
superexponential then $d$-regularity  is also necessary for the family $(
a_\alpha^{-d})_{\alpha\in\ell_A}$ to be summable.
\end{Pn}

\begin{proof}[\indent Proof]
We first prove the theorem for the case $d=1$.
It follows from \eqref{eq:ineq} that for every $\alpha \in
\ell_A$,
\[
\prod_{n \in A} a_{e_n}^{\alpha_n} \leq a_\alpha.
\]
Therefore,
\[
\begin{split}
\sum _{\alpha \in \ell_A} a_\alpha^{-1}
&\leq \sum _{\alpha \in \ell_A}\prod_{n \in A}a_{e_n}^{-\alpha_n} \\
&=\prod_{n \in A} \sum _{\alpha_n=0}^{\infty}a_{e_n}^{-\alpha_n}\\
&=\prod_{n \in A} \frac 1{1-a_{e_n}^{-1}} \\
&=\prod_{n \in A}\left(1+\frac 1{a_{e_n}-1}\right) < \infty
\end{split}
\]
If $a_\alpha a_\beta = a_{\alpha +\beta}$, then $\forall \alpha
\in \ell_A$, $\prod_{n \in A} a_{e_n}^{\alpha_n} = a_\alpha $.
Therefore,
\[
\sum _{\alpha \in A} a_\alpha^{-1} = \prod_{n \in A} \left(1+\frac
1{a_{e_n}-1}\right),
\]
which converges if and only if  $\sum _{n \in A} \frac
1{a_{e_n}-1}< \infty$. In case $d>1$, we take $a^d$ instead
of $a$.
\end{proof}

When $A=\mathbb N$ and $a_\alpha$ is given by \eqref{eq:zhang} we
obtain as a corollary a result of Zhang, proved in 1992, see
\cite{zhang}, \cite{MR1408433}. Zhang's proof uses Abel's convergence
test. The result
itself is necessary in order to present some important properties
of Kondratiev spaces of stochastic test functions and stochastic
distributions.

\begin{Cy}[Zhang \cite{zhang}] Let $d\in\mathbb N$.
$\sum _{\alpha \in \ell} (2  \mathbb{N})^{-d \alpha} < \infty$ if
and only if $d>1$.
\end{Cy}

\begin{proof}[\indent Proof]
We take $\ell_A=\ell$. Thus
\[
\sum _{\alpha \in \ell} (2\mathbb{N})^{-d \alpha} <
\infty\quad\iff \quad\sum _{n \in \mathbb{N}} \frac 1{(2n)^d-1}<
\infty,
\]
which is true if and only if $d>1$.
\end{proof}


The spaces $\mathscr S^\prime$ and $\mathcal S_{-1}$ are strong
dual of Fr\'echet spaces. Namely, $\mathscr S^\prime$ is the
strong dual of the Schwartz space $\mathscr S$ of Hermite series
\begin{equation}
\label{eq:hermite}
f(x)=\sum_{n=0}^\infty f_n\xi_n(x),
\end{equation}
where $(\xi_n)_{n\in\mathbb N_0}$ denote the Hermite functions, and
the coefficients $f_n$ are complex numbers, and such that
\[
\sum_{n \in \mathbb{N}_0} {
 (n+1)^{2p}|f_n| ^2<\infty} \text{ for
all $p \in \mathbb{N}$,}
\]
and $\mathcal S_{-1}$ is the dual of the Kondratiev space
$\mathcal S_{1}$ of stochastic test functions which can be seen
as families of complex numbers $(f_\alpha)_{\alpha\in\ell}$
indexed by $\ell$ and such that
\[
\sum_{\alpha \in \ell} {(\alpha!)^2(2
\mathbb{N})^{\alpha p}|f_\alpha| ^2 <\infty} \text{ for all  $p \in
\mathbb{N}$},
\]
where
\begin{equation}\label{factorial}
(2\mathbb N)^\alpha=\prod_{k=1}^\infty (2k)^{\alpha_k} \quad \text{and}\quad  \alpha!=\prod_{k=1}^\infty (\alpha_k !).
\end{equation}
The triples $(\mathscr S,\mathbf L_2(\mathbb R,dx),\mathscr
S^\prime)$ and $(\mathcal S_1,\mathcal W,\mathcal S_{-1})$, where
$\mathcal W$ denotes the white noise space, are Gelfand triples.
We now define a wide family of nuclear topological vector spaces,
which includes  $\mathscr S^\prime$ and $\mathcal S_{-1}$, and
which is closed under tensor products, as dual of certain
Fr\'echet spaces. \\

Let $a:\ell_A \to \mathbb R$ be a positive function, that is $\alpha \mapsto a_\alpha$ where $a_\alpha>0$ for any $\alpha  \in \ell_A$. We denote the weighted Hilbert space with respect to $a$ by
\begin{equation}\label{Hil_def}
\ell^2_a=\left\{(\varphi_\alpha)_{\alpha \in \ell_A} :
\sum_{\alpha \in \ell_A} |\varphi_\alpha|^2
a_\alpha<\infty \right\}.
\end{equation}
When $a$ is the constant function $1$ we simply denote $\ell^2=\ell^2_1$ (i.e $\ell^2=\ell^2(\ell_A)$).
We define
the countably normed space
\begin{equation}
\label{eq:Fdef}
\mathcal F_{a}=\left\{(\varphi_\alpha)_{\alpha \in \ell_A} :
\sum_{\alpha \in \ell_A} |\varphi_\alpha|^2
a_\alpha^{p}<\infty \text{ for all } p \in \mathbb{N} \right\} = \bigcap_{p \in \mathbb N} \ell_{a^p}.
\end{equation}

\begin{Tm}\label{Tm:nuc}
Let $a:\ell_A \to \mathbb R$ be a positive function.
\begin{enumerate}[(a)]
\item The space $\mathcal F_{a}$ endowed with the topology
defined by the norms
\[
\|\varphi\|_p^2= \sum_{\alpha
\in \ell_A} |\varphi_\alpha|^2 a_\alpha^{p},\quad
p=1,2,\ldots
\]
is a Fr\'echet space.
\item If $a>1$ (that is, for any $\alpha \in \ell_A$, $a_\alpha>1$), then the space $\mathcal F_{a}$ is continuously included in $\ell^2$.\\
\item The space $\mathcal F_{a}$ is
nuclear if and only if there exists $d \in \mathbb N$ such that
\begin{equation}\label{eq:nuc}
\sum _{\alpha \in \ell_A} a_\alpha^{-d} < \infty
\end{equation}
\end{enumerate}
\end{Tm}

\begin{proof}[\indent Proof]
\mbox{}\\
\begin{enumerate}[(a)]
\item
By definition, $\mathcal F_{a}=\bigcap_{p\in \mathbb N} \ell_{a^p}$, and thus it is a Fr\'echet space.
We first note that the norms are non decreasing and compatible
(i.e. every sequence which is a Cauchy sequence with respect to the two
norms and converges to zero with respect to one of them,
converges to
zero also with respect to the second). See \cite[p. 17-18]{GS2_english}\\
\item We have
\[
\|\varphi\|_{\ell^2}^2= \sum_{\alpha
\in \ell_A}|\varphi_\alpha|^2  \leq \sum_{\alpha \in
\ell_A}|\varphi_\alpha|^2 a_\alpha^p =
\|\varphi\|_{\ell_{a^p}}^2.
\]

\item Defining $\delta_\alpha = (\delta_{\alpha,\beta})_{\beta \in
\ell_A}$ such that $\delta_{\alpha, \beta}=0$ if $\alpha \neq
\beta$  and $\delta_{\alpha, \beta}=1$ if $\alpha=\beta$, it is
clear that $(\delta_\alpha a_\alpha^{-p/2} b_\alpha)_{\alpha \in
\ell_A}$ is an orthonormal base of $\ell_{a^p}$.
Now, $\mathcal F_{a}=\bigcap_{p \in \mathbb{N}}
\ell_{a^p}$ is nuclear if and only if for every
$p \in \mathbb N$ there exists $q>p$ such that the natural
embedding $\iota : \ell_{a^q} \to \ell_{a^p}$
is Hilbert-Schmidt. The equality
\[
\begin{split}
{\rm tr}~(\iota ^* \iota) &=\sum_{\alpha \in \ell_A}\langle \iota
^* \iota(\delta_\alpha a_\alpha^{-q/2}b_\alpha),
(\delta_\alpha a_\alpha^{-q/2}b_\alpha) \rangle _q\\
&=\sum_{\alpha \in \ell_A}\|
\iota(\delta_\alpha a_\alpha^{-q/2}b_\alpha)\|_p ^2 \\
&=\sum_{\alpha \in \ell_A}a_\alpha^{-(q-p)}\| \iota(\delta_\alpha
a_\alpha^{-p/2}b_\alpha)\|_p ^2=\sum_{\alpha \in
\ell_A}a_\alpha^{-(q-p)}
\end{split}
\]
yields that $\mathcal F_a$ is nuclear if and only if for any $p$ there exists $q>p$ such that
\begin{equation} \label{eq:nuc2}
\sum_{\alpha \in \ell_A}a_\alpha^{-(q-p)}<\infty.
\end{equation}
If \eqref{eq:nuc} holds for some $d \in \mathbb N$, then setting
$q=p+d$ leads to the requested result. On the opposite direction,
if $\mathcal F_a$ is nuclear, then for $p=1$ there exists $q>p$
such that \eqref{eq:nuc2} holds. Setting $d=q-p=q-1$ yields the requested result.
\end{enumerate}
\end{proof}


\begin{Dn}
\label{paris:distribution_space}
Let $a$ be a $d$-regular admissible positive function. The space $\mathcal F_{a}^\prime$ is called a $d$-regular 
admissible space. The space is called regular admissible if it is $d$-regular admissible for
some $d\in\mathbb N$.
\end{Dn}

We denote by $(\ell^2_{a^p})^\prime$ the dual of $\ell^2_{a^p}$.
Then
\[
(\ell^2_{a^1})^\prime \subseteq (\ell^2_{a^2})^\prime \subseteq
\cdots \subseteq (\ell^2_{a^p})^\prime \subseteq \cdots \subseteq
\mathcal F_{a}^\prime,
\]
and the dual space $\mathcal F_{a}^\prime$ is the union of the
increasing sequence of the spaces $(\ell^2_{a^p})^\prime$, i.e.,
\[
\mathcal F_{a}^\prime=\bigcup_{p \in \mathbb{N}}
(\ell^2_{a^p})^\prime.
\]
See \cite[p. 35-36]{GS2_english}. Since a Fr\'{e}chet space is nuclear if
and only if its strong dual is nuclear, $\mathcal F_{a}^\prime$
is also nuclear.

\begin{Pn}
$\mathcal F_{a}^\prime$ can be viewed as
\[
\left\{(f_\alpha)_{\alpha \in \ell_A} :
\sum_{\alpha \in \ell_A}
|f_\alpha| ^2 a_\alpha^{-p}<\infty \text{ for some }p \in
\mathbb{N} \right\}.
\]
\end{Pn}

\begin{proof}[\indent Proof]
Let $f \in (\ell^2_{a^p})^\prime$. It follows from Riesz's
representation theorem that there exists $\psi=(\psi_\alpha) \in \mathscr \ell^2_{a^p}$ with
$\|\psi\|_{\ell^2_{a^p}}=\|f\|_{(\ell^2_{a^p})^\prime}$ and such that
\[
f(\cdot)=\langle \cdot,\psi\rangle_{\ell^2_{a^p}}.
\]
Thus, for any $\varphi=(\varphi_\alpha) \in \ell^2_{a^p}$
\[
f(\varphi)=\langle \psi,\varphi \rangle_{\ell^2_{a^p}}
=\sum_{\alpha \in \ell_A} \varphi_\alpha \overline{\psi_\alpha}
a_{\alpha}^p.
\]
Setting $f_\alpha= \psi_\alpha a_{\alpha}^p$, we have
\[
\sum_{\alpha \in \ell_A}|f_\alpha|^2a_{\alpha}^{-p}
=\sum_{\alpha \in \ell_A}|\psi_\alpha|^2 a_{\alpha}^p = \| \psi
\|^2_{\ell^2_{a^p}}=\|f\|^2_{(\ell^2_{a^p})^\prime}.
\]
Moreover, for any $(f_\alpha)$ subjects to $\sum_{\alpha \in \ell_A}
|f_\alpha|^2a_{\alpha}^{-p}<\infty$
\[
(f_\alpha) \mapsto \left( (\varphi_\alpha) \mapsto \sum_{\alpha \in \ell_A}
\varphi_\alpha \overline{f_\alpha} \right)
\]
maps $(f_\alpha)$ to a continuous linear functional over $\ell^2_{a^p}$,
and any composition of this mapping with $f \mapsto (f_\alpha)$, which was
described before, yields the appropriate identity. Hence,
\begin{equation}
\label{norm}
(\ell^2_{a^p})^\prime \cong \ell^2_{a^{-p}}=
\left\{(f_\alpha)_{\alpha \in
\ell_A}:\, \sum_{\alpha \in
\ell_A}|f_\alpha|^2a_{\alpha}^{-p}<\infty \right\}.
\end{equation}
Thus, $\mathcal F_{a}^\prime$ can be viewed as
\[
\bigcup_{p \in \mathbb N} \ell^2_{a^{-p}}=\left\{(f_\alpha)_{\alpha \in \ell_A} :
\sum_{\alpha \in \ell_A}
|f_\alpha| ^2 a_\alpha^{-p}<\infty \text{ for some }p \in
\mathbb{N} \right\}.
\]
\end{proof}

We note that the inner product $\langle
\cdot,\cdot \rangle _{\ell^2}$ coincides with the antilinear duality
$\langle \cdot,\cdot
\rangle _{\mathcal
F_{a}^\prime,\mathcal F_{a}}$, whenever it makes sense.
Considering the inclusion of dual spaces $(\ell^2)^\prime$ in
$\mathcal F_{a}^\prime$, using Riesz theorem we have that,
\[
\mathcal F_{a} \subseteq \ell^2 \subseteq \mathcal
F_{a}^\prime.
\]
It is clear that the second inclusion is also continuous, since
\[
\|f\|_{\mathscr
(\ell^2_{a^p})^\prime}^2=\sum_{\alpha \in \ell_A}|f_\alpha|^2
a_\alpha^{-p} \leq \sum_{\alpha \in \ell_A}|f_\alpha|^2
 = \|f\|_{\ell^2}^2.
\]

\begin{Pn}
$(\mathcal F_{a}, \ell^2(\ell_A), \mathcal F_{a}^\prime)$ is a
Gelfand triple.
\end{Pn}


\section{V{\aa}ge spaces - a new family of topological $\mathbb C$-algebras}
\label{sec:3} \setcounter{equation}{0}
Considering the monoid $\ell_A$ for $A \subseteq \mathbb N$,
it has the property that for any $\gamma \in\ell_A$
\[
\{(\alpha,\beta)\in \ell_A^2 : \alpha+\beta=\gamma\} \text { is finite}.
\]
Thus it gives rise to the total $\mathbb C$-algebra of the monoid $\ell_A$, namely $(\mathbb C^{\ell_A},+,*)$, where $*$ denotes the convolution multiplication defined by
\begin{equation}\label{eq:conv}
f*g
=\left(\sum_{\alpha+\beta=\gamma}f_\alpha g_{\beta}\right)_{\gamma \in \ell_A} \text{ for all } f,g \in \mathbb C^{\ell_A}.
\end{equation}
We usually omit the $*$, and simply write $fg$ instead of $f*g$.
This algebra is an integral domain and a topological ring with respect to the product topology of $\mathbb C^\ell_A$.
Denoting $x_n=(\delta_{e_n,\beta})_{\beta \in \ell_A} \in \mathbb C^{\ell_A}$, $x=(x_n)_{n \in A}$ and $x^\alpha=\prod_{n \in A}x_n^{\alpha_n}$, we obtain $x^\alpha=(\delta_{\alpha,\beta})_{\beta \in \ell_A}=\delta_\alpha$. Therefore
\[
(f_\alpha)_{\alpha \in \ell_A}= \sum_{\alpha \in \ell_A}f_\alpha x^\alpha.
\]
Thus, we obtain $\mathbb C^{\ell_A}=\mathbb C[[x]]$ and the convolution can be considered simply as formal series multiplication.
For more details about total algebras of monoids see Bourbaki \cite[p. 454]{Bourbaki_Algebra}.



The proof of the following proposition is straightforward and will
be omitted.
\begin{Pn} Let $f,g$ and $h$ be in $\mathbb C^{\ell_A}$. Then, it
holds that:
\begin{enumerate}[(a)]
\item $f  g  = g  f $.
\item $(f  g) h = f ( g h)$.
\item $ f ( g + h) =f  g + f  h$.
\end{enumerate}
\end{Pn}

Recall that we introduced regular admissible spaces in Definition
\ref{paris:distribution_space}. When it does not lead to any
confusion, we simply denote by $\| \cdot \|_p$ instead of $\|
\cdot \|_{(\ell^2_{a^{p}})'}=\| \cdot \|_{\ell^2_{a^{-p}}}$.

\begin{Dn}
A regular admissible space $\mathcal F_{a}^\prime=\bigcup_{p=1}^\infty
\ell^2_{a^{-p}}$ is called a V{\aa}ge space if there is $e \in
\mathbb{N}$ such that for every $p\in\mathbb N$ and for every
$p \geq q+e$
\[
\|f g\|_{p} \leq A(p-q)\|f\|_{q}\|g\|_{p}
\]
for all $f \in \mathscr \ell^2_{a^{-p}}$ and $g \in \ell^2_{a^{-p}}$,
where $A(p-q)$ is a finite positive number.
If $\mathcal F_{a}'$ is a  V{\aa}ge space, we call the minimal $e$ with
this property the {\sl index} of the space.
\end{Dn}

We note that in particular, a V\aa ge space $\mathcal F_a'$
is a subalgebra of $(\mathbb C^{\ell_A},+,*)$, i.e.,
closed under convolution (and clearly under addition).

\begin{Tm}\label{Tm:VageSpace}
A $d$-regular admissible space $\mathcal F_{a}^\prime$ is a V{\aa}ge space if and
only if
\[
a_\alpha a_\beta \leq a_{\alpha +\beta},\quad\forall\alpha, \beta
\in \ell_A,
\]
i.e., $a:\ell_A \to \mathbb R$ is superexponential. Its
index is then less or equal to $d$.
\end{Tm}

\begin{proof}[\indent Proof]
We follow the argument in \cite[p. 118]{MR1408433}.
First we assume that for all $\alpha, \beta \in
\ell_A$, $a_\alpha a_\beta \leq a_{\alpha +\beta}$.
Since $\mathcal F_{a}'$ is a regular admissible space, $a_0=1$, $a_{e_n}>1$ for
all $n \in A$, and  $\sum _{n \in A} \frac 1{a_{e_n}^d-1}< \infty$
for some $d \in \mathbb{N}$. Therefore by Proposition \ref{Pn:Converges} we obtain for any $p-q\geq d$
$\sum_{\alpha \in \ell_A} a_{\alpha}^{-(p-q)} \leq \sum_{\alpha \in
\ell_A} a_{\alpha}^{-d}<\infty$.
We denote
\[
A(p-q)=\left(\sum_{\alpha \in \ell_A} a_\alpha^{-(p-q)}\right)^{\frac 12}.
\]
Now, supposed that $f \in \ell^2_{a^{-q}}$
and $g\in \ell^2_{a^{-p}}$, for some $p \geq q+d$.
Then,
\[
\begin{split}
\|f g\|_{p} ^2
& =\sum_{\gamma \in \ell_A}\left|\sum_{\alpha \leq \gamma} f_\alpha
g_{\gamma-\alpha}a_\gamma^{-p/2}\right|^2\\
& \leq \sum_{\gamma \in \ell_A}\left(\sum_{\alpha \leq \gamma}|
f_\alpha| a_{\alpha}^{-p/2} |g_{\gamma-\alpha}|
a_{\gamma-\alpha}^{-p/2}\right)^2\\
& \leq \sum_{\gamma \in \ell_A}\left(\sum_{\alpha,\alpha' \leq \gamma}|
f_\alpha|a_{\alpha}^{-p/2}|f_{\alpha'}| a_{\alpha'}^{-p/2}
|g_{\gamma-\alpha}|a_{\gamma-\alpha}^{-p/2}|g_{\gamma-\alpha'}
|a_{\gamma-\alpha'}^{-p/2}\right)\\
& \leq \sum_{\alpha,\alpha' \in \ell_A}\left(| f_\alpha|a_{\alpha}^{-p/2}
|f_{\alpha'}| a_{\alpha'}^{-p/2} \sum_{\gamma \geq \alpha,\alpha'}
|g_{\gamma-\alpha}|a_{\gamma-\alpha}^{-p/2}|g_{\gamma-\alpha'}
|a_{\gamma-\alpha'}^{-p/2} \right)\\
& \leq \left(\sum_{\beta \in \ell_A}|f_\beta| a_{\beta}^{-p/2}
\right)^2 \left(\sum_{\beta \in \ell_A}|g_{\beta}|^2 a_{\beta}^{-p}\right)^{\frac 12}
 \left(\sum_{\beta \in \ell_A}|g_{\beta}|^2 a_{\beta}^{-p}\right)^{\frac 12}\\
&\leq \left(\sum_{\beta \in \ell_A} a_{\beta}^{-(p-q)}\right)
\left(\sum_{\beta \in \ell_A} |f_\beta|^2 a_{\beta}^{-q}\right) \left(\sum_{\beta
 \in \ell_A}|g_{\beta}|^2 a_{\beta}^{-p}\right)\\
&=\left(A(p-q)\right)^2\|f\|_{q}^2\|g\|_{p}^2.
\end{split}
\]
Thus, $\mathcal F_{a}$ is a V\aa ge space of index which is less than or equal to $d$.
On the other direction, assuming that $\mathcal F_{a}'$ is a V{\aa}ge space of
index $e$.
Let $q$ be a natural number and $p=q+e$. Then,
\[
\|x^{\alpha+\beta}\|_{p}=\|x^\alpha
x^\beta\|_{p} \leq
A(p-q)\|x^\alpha\|_{q}\|x^\beta\|_{p}.
\]
Therefore,
$a_{\alpha+\beta}^{-p} \leq A(e) a_\alpha^{-p} a_\beta^{-q}$.
Then,  $ A(e)^{-1/p} a_\alpha a_\beta^{q/p} \leq
a_{\alpha+\beta}$. Thus, we obtain $ a_\alpha a_\beta \leq a_{\alpha+\beta}$
as $p$ goes to infinity.
\end{proof}

\begin{Cy}
A V{\aa}ge space is nuclear.
\end{Cy}
\begin{proof}[\indent Proof]
If $\mathcal F_{a}^\prime$ is a V{\aa}ge space then
$\sum_{\alpha \in \ell_A} a_\alpha^{-d}<\infty$. Applying
Theorem \ref{Tm:nuc}, $\mathcal F_{a}$ is nuclear.
However, since it is a  Fr\'echet space, $\mathcal
F_{a}^\prime$ is also nuclear.
\end{proof}

We now show that the convolution is continuous in the strong topology of a V\aa ge space. Before that, we need the following proposition.

\begin{Pn}\label{Pn:net}
Let $\mathcal F_a'$ be a regular admissible space and let $(f_\lambda)$ be a net in $\mathcal F_{a}'$. Then  $f_\lambda
\to f$ in the strong topology if and only if  there exists $p
\in \mathbb{N}$ such that $f_\lambda,f \in \ell^2_{a^{-p}}$ and $f_\lambda
\to f$ in the strong topology of $\mathscr \ell^2_{a^{-p}}$.
\end{Pn}
\begin{proof}[\indent Proof]
Suppose $f_\lambda \to f$ in the strong topology of $\mathcal
F_{a}'$. In particular, $\{f_\lambda\}_{\lambda \in \Lambda
}\cup \{f\}$ is strongly bounded. Therefore, there exists $p \in
\mathbb{N}$ such that
\[
\{f_\lambda\}_{\lambda \in \Lambda}\cup
\{f\} \subseteq \ell^2_{a^{-p}}
\]
(see \cite[\S 5.3 p.
45]{GS2_english}).
Let $B$ be a bounded set in $\ell^2_{a^{p}}$, then $B \cap
\mathcal F_{a}$
  is a dense subset of $B$. Therefore,
\[
\sup_{\varphi \in B}
  |f_\lambda(\varphi) - f(\varphi)|=\sup_{\varphi \in B \cap
   \mathcal F_{a}}|f_\lambda(\varphi) - f(\varphi)| \to 0.
\]
Thus,
$f_\lambda \to f$ in the strong topology of $\ell^2_{a^{-p}}$.
 The opposite direction is clear.
\end{proof}

\begin{Tm}
Let $\mathcal F_{a}'$ be a V{\aa}ge space. Then the convolution
is a continuous function $\mathcal F_{a}' \times \mathcal F_{a}' \to
\mathcal F_{a}'$ in the strong topology. Hence $(\mathcal F_{a}',+,*)$ is a topological $\mathbb C$-algebra.
\end{Tm}
\begin{proof}[\indent Proof]
Assuming $((f_\lambda,g_\lambda))_{\lambda \in \Lambda}$ is a net which converges
to $(f,g)$ in the strong topology of $\mathcal F_{a}' \times \mathcal F_{a}'$, then in
particular, $f_\lambda \to f$ and $g_\lambda \to g$ in the strong
topology of $\mathcal F_{a}'$. According to Proposition \ref{Pn:net}, there
exist $p,q \in \mathbb{N}$ such that $f_\lambda,f \in \ell^2_{a^{-q}}$
and $g_\lambda,g \in \ell^2_{a^{-q}}$ where $f_\lambda \to f$ in the
strong topology of $\ell^2_{a^{-p}}$ and $g_\lambda \to g$ in the strong
topology of $\ell^2_{a^{-q}}$. We may assume that $p \geq q+d$. Since
$\mathcal F_{a}'$ is a V{\aa}ge space, $f g_\lambda=f*g_\lambda, f g=f*g
\in \ell^2_{a^{-p}}$, and $*: \ell^2_{a^{-q}} \times \ell^2_{a^{-p}}
 \to
\mathscr \ell^2_{a^{-p}}$
is continuous.\\
Since $(f_\lambda,g_\lambda) \to (f,g)$ in the
strong topology of $\ell^2_{a^{-q}} \times \ell^2_{a^{-p}}$, $f_\lambda
g_\lambda \to f g$ in the strong topology of $\ell^2_{a^{-p}}$.
Again, using Proposition \ref{Pn:net}, we have that $f_\lambda
 g_\lambda \to f  g$ in the strong topology of
$\mathcal F_{a}'$. Thus the convolution is strongly continuous.
\end{proof}

We do not know if the convolution is continuous in the weak topology.
We end this section with the
weak topology analogue of Proposition \ref{Pn:net}.

\begin{Pn}
Let $(f_\lambda)$ be a net in $\mathcal F_{a}'$. Then
$f_\lambda \to f$ in the weak topology if and only if
there exists $p \in \mathbb{N}$ such that $f_\lambda,f \in \ell^2_{a^{-p}}$
and $f_\lambda \to f$ in the weak topology of $\ell^2_{a^{-p}}$.
\end{Pn}
\begin{proof}[\indent Proof]
Suppose $f_\lambda \to f$ in the weak topology .In particular,
$\{f_\lambda\}\cup \{f\}$ is weakly bounded, and thus strongly bounded
(see \cite[p. 48]{GS2_english}). Therefore,
there exists $p \in \mathbb{N}$ such that $\{f_\lambda\}
\cup \{f\} \subseteq \mathscr \ell^2_{a^{-p}}$. Moreover, $f_\lambda \to f$
 pointwise on a dense subset of $\ell^2_{a^{p}}$, that is $\mathcal F_{a}$.
 Let $\epsilon>0$ and $\varphi \in \ell^2_{a^{p}}$. We may choose
 $\psi \in \mathcal F_{a}$ such that $\|\varphi - \psi \|_p<
 \frac{\epsilon}{2(\|f\|+\sup_\lambda\|f_\lambda\|)}$,
 and $\lambda_0 \in \Lambda$ such that for all $\lambda
  \geq \lambda_0$, $|f_\lambda(\psi)-f(\psi)|<\frac{\epsilon}{2}$. Therefore,
\[
\begin{split}
|f_\lambda(\varphi)-f(\varphi)|
&\leq |f_\lambda(\varphi)-f_\lambda(\psi)|+|f_\lambda(\psi)
-f(\psi)|+|f(\psi)-f(\varphi)| \\
&\leq (\|f_\lambda\|+\|f\|)\|\varphi-\psi\|+|f_\lambda(\psi)-f(\psi)| \\
&<\frac{\epsilon}{2}+\frac{\epsilon}{2}=\epsilon.
\end{split}
\]
Thus $f_\lambda \to f$ in the weak topology of $\ell^2_{a^{-p}}$.
The opposite direction is clear.
\end{proof}

\section{Invertible elements and power series}
\label{sec:4}
\setcounter{equation}{0}
The $\mathbb C$-algebra $(\mathcal F_{a}',+,*)$ is in particular a (unit) ring. We denote it by $\bigR$.
\begin{Dn}
Let $f =(f_\alpha)_{\alpha \in \ell_A} \in \bigR$.
Then, $f_0 \in \mathbb{C}$ is called the generalized
expectation of $f$ and is denoted by $E[f]$.
\end {Dn}
From this definition we have that
\[
E[f  g] = E[f]E[g], \quad \forall f,g \in \bigR.
\]
We note that $E: \bigR \to \mathbb{C}$ is a homomorphism which
maps $1_\bigR$ to $1_\mathbb{C}$. In the sequel, we will see it is
the only homomorphism with this property (see Proposition \ref{prop: homo}).

\begin{Pn} \label{Pn:limF}
Let $M$ be a positive  number. Then, for any $f \in \bigR$ such
that $E[f]=0$ there is $q \in \mathbb{N}$ such that $\|f\|_{q}<M$.
\end{Pn}

\begin{proof}[\indent Proof]
Let $f=(f_\alpha) \in \ell^2_{a^{-p}}$ with $f_0=0$.
Since $a_{e_n}>1$ for $n=1,2,\ldots$, we have
\[
a_\alpha=\prod_{n \in A}a_{e_n}^{\alpha_n}>1 \quad \text{ for all } \alpha
\neq (0,0,\ldots).
\]
Therefore,
for all $\alpha \in \ell_A$ $\lim_{q \to \infty}|f_\alpha|^2a_\alpha^{-q}=0$
(recall $f_0=0$)
and for all $q>p$, $|f_\alpha|^2a_\alpha^{-q} \leq |f_\alpha|^2a_\alpha^{-p}$,
whereas $\sum_{\alpha \in \ell_A}|f_\alpha|^2a_\alpha^{-p}=\|f\|_p^2<\infty$.
Thus, the dominated convergence theorem implies
\[
\lim_{q \to \infty}\|f\|_q^2
=\lim_{q \to \infty} \sum_{\alpha\in\ell_A} |f_\alpha|^2 a_\alpha^{-q}
=\sum_{\alpha\in\ell_A}\lim_{q \to \infty}  |f_\alpha|^2 a_\alpha^{-q}
=0.
\]
\end{proof}

\begin{Dn}
The $n$-th convolution power $f^n=f^{*n}$ of $f$ is defined
inductively as follows:
$$
f^{n}=
  \begin{cases}
   1 & \text{if } n=0\\
   f f^{(n-1)} & \text{if } n>0
  \end{cases}
$$
\end{Dn}

\begin{Pn} \label{Pn:power}
Let $\bigR$ be a V{\aa}ge space of index $d$ and
let $f$ be in $\ell^2_{a^{-p}}$. Then $\forall n \in \mathbb{N}$,
$f^{n} \in \ell^2_{a^{-(p+d)}}$. Moreover, $\|f^n\|_{p+d}
\leq A(d)^n  \|f\|_{p}^n$.
\end{Pn}

\begin{proof}[\indent Proof]
Obviously, $f^0=1 \in \ell^2_{a^{-(p+d)}}$, and
$\|f^{ 0}\|_{p+d} = A(d)^0  \|f\|_{p}^0$.\\
By induction, if we assume $f^{n} \in \bigR$, we get
$f^{(n+1)}=ff^{n} \in \bigR$, and
\[
\begin{split}
\|f^{(n+1)}\|_{p+d}
&=\|f f^{n}\|_{p+d}\\
&\leq A(d) \|f\|_{p} \|f^{ n}\|_{p+d} \\
&\leq A(d)^n \|f\|_{p}^{n+1} <\infty
\end{split}
\]
\end{proof}

More generally, given a polynomial
$p(z)=\sum_{n=0}^N p_n z^n$  ($p_n \in \mathbb{C}$), we define its convolution version
$ p:\bigR \rightarrow \bigR$
by
\[
p(f)=\sum_{n=0}^N p_n f^{n}
\]
By Proposition \ref{Pn:power}, we have that $p(f) \in \bigR$ for $f \in
\bigR$. The following proposition considers the case of power
series.

\begin{Pn}\label{Pn:series}
Let $\phi(z)=\sum_{n \in \mathbb{N}} \phi_n z^n$ be a power
series (with complex coefficients)
which converges absolutely in the open disk with radius
$R$. Then for any $f \in \bigR$ such that $|E[f]|<\frac{R}{A(d)}$  it holds that
\[
\phi(f)=\sum_{n \in \mathbb{N}} \phi_n f^{n}
\in \bigR.
\]
\end{Pn}

\begin{proof}[\indent Proof]
Applying Proposition \ref{Pn:limF}, there exists $q$ such that $\|f-E(f)\|_q<\frac{R}{A(d)}-|E[f]|$. Therefore,
\[
\|f\|_q \leq \|f-E(f)1_{\bigR}\|_q+|E(f)|<\frac{R}{A(d)}.
\]
Then by Proposition \ref{Pn:power}, for all $p \geq q+d$,
\[
\begin{split}
\sum_{n \in \mathbb{N}} |\phi_n|\| f^{n}\|_{p}
&\leq \sum_{n \in \mathbb{N}} |\phi_n|A(d) ^n \|f\|_{q}^n\\
&= \sum_{n \in \mathbb{N}} |\phi_n|(A(d) \|f\|_{q})^n\\
&<\infty.
\end{split}
\]
Since $\ell^2_{a^{-p}}$ is a Hilbert space, $\phi(f)=\sum_{n
\in \mathbb{N}} \phi_n f^{n} \in \ell^2_{a^{-p}}$. Thus,
$\phi(f) \in \bigR$.
\end{proof}

\begin{Pn} \label{prop: inv}
The element $f \in \bigR$ is invertible  if and only if $E[f]$ is invertible.
\end{Pn}

\begin{proof}[\indent Proof]
If $E[f] \neq 0$, we can assume that $E[f] =1$. By the Proposition \ref{Pn:series}
we have that $\sum_{n \in \mathbb{N}} (1-f)^{
n} \in \bigR$. Furthermore,
\[
f\left(\sum_{n \in \mathbb{N}} (1-f)^{
n}\right)=1.
\]
Conversely, assume $f$ invertible. Then there exists $f^{-1} \in
\bigR$ such that $f  f^{-1} =1$. Hence,
$E[f]E[f^{-1}]=E[f  f^{-1}]=1$.
\end{proof}

\begin{Pn}\label{prop: homo}
Let $\bigR$ be a V{\aa}ge space. Then the following properties hold:
\begin{enumerate}[(a)]
\item $GL(\bigR)$ is open.
\item The spectrum of $f \in \bigR$, $\sigma(f)=
\{\lambda \in \mathbb C : f-\lambda 1_\bigR \text{ is not invertible }\}$
 is the singleton $\{E[f]\}$.
\item $E$ is unique as a homomorphism $\bigR \to \mathbb{C}$
mapping $1_\bigR$ to $1_\mathbb{C}$.
\end{enumerate}
\end{Pn}
\begin{proof}[\indent Proof]
$\quad$\\
\begin{enumerate}[(a)]
\item
By Proposition \ref{prop: inv}, we have  that $\{f \in \bigR : E[f] \neq
0\}$ is the set of all invertible elements in $\bigR$. In other
words, $GL(\bigR)=E^{-1}(GL(\mathbb{C}))$. In particular, since
$E$ is continuous, $GL(\bigR)$ is open.
\item
Clearly, $f-\lambda1_\bigR$ does not have an inverse
if and only if $\lambda=E(f)$.
\item
Let $\varphi:\bigR \to \mathbb C$ be a  homomorphism
mapping $1_\bigR$ to $1_\mathbb{C}$ and let $f \in \bigR$.
Since $\varphi \left(f-\varphi(f) 1_\bigR   \right)=0$,
$\varphi(f) \in \sigma(f)$, that is $\varphi(f)=E[f]$.
\end{enumerate}
\end{proof}

The notion of rational functions plays a key role in the theory of
linear systems over commutative rings. See Section \ref{sec:8}
below for a discussion and references. Here, using Proposition
\ref{prop: inv} we define rational functions with coefficients in
$\bigR$ as in \cite{alp} and \cite[Section 3] {aa_goh} as
elements in the quotient ring of $\mathcal R[z]$. Other
characterizations can be given in terms of finite dimensional
backward shift invariant spaces and realizations. See \cite{alp}.
We mention that in the non-commutative case similar characterisations
occur. See for instance \cite{acs1,MR2124899}.

\begin{Dn}
\label{def4}
A rational function with coefficients in $\bigR^{n\times m}$ is an
expression of the form
\begin{equation}
\label{eq:rat1}
R(z)= p(z) ( q(z))^{-1}
\end{equation}
where $ p\in(\bigR[z])^{n\times m}$, and $
0\neq q\in\bigR[z]$.
\end{Dn}

Note that, for any $f\in\mathcal R$ such that $E(q(f))\not =0$, $R(f)$ is well defined
as an element of $\mathcal R^{n\times m}$.\\

In Section \ref{sec:8} we discuss some problems in the theory of
linear systems using rational functions with coefficients in $\bigR$.\\

To conclude this section, we remark that one may consider
classical interpolation problems (of the kind developed in \cite{bgr}) in
the present setting. For instance one can consider the following
Nevanlinna-Pick interpolation problem:

\begin{Pb}
\label{pb:NP}
Given $N\in\mathbb N$ and $N$ pairs of points $(a_1,b_1),\ldots, (a_N,b_N)$ in
$\mathcal R^2$, find all power series  $\phi$ such that
\[
\phi(a_j)=b_j,\quad j=1,\ldots, N,
\]
and such that, moreover, the function $z\mapsto E(\phi(z))$ is analytic and contractive
in the open unit disk.
\end{Pb}

This problem, as well as more general interpolation problems, have been studied
in \cite{aa_goh} when $\mathcal R$ is the Kondratiev
space $\mathcal S_{-1}$. Two important tools used there, and which are
still valid in the setting of V\aa ge spaces are the permanence of
algebraic identities (see
\cite[p. 456]{MR1129886}) and the definition and properties of analytic
functions with
values in the dual of a countably normed Hilbert space.


\section{Tensor product of V{\aa}ge spaces}
\setcounter{equation}{0} \label{sec:5} When one considers the
tensor product $E \otimes F$ of two locally convex Hausdorff
spaces $E,F$, some "natural" topologies may be considered. Two
such topologies are the $\pi$-topology and the $\epsilon$-topology
(see \cite{Groth55}, \cite[Chapter 43, p. 434]{Treves67}). These
topologized tensor products of $E$ and $F$ are denoted
respectively by $E \otimes_\pi F$ and $E \otimes_\epsilon F$. The
completions of the tensor product of $E$ and $F$ with respect to
 the $\pi$-topology and the $\epsilon$-topology, will then be denoted
 by $E \hat \otimes_\pi F$ and $E \hat \otimes_\epsilon F$ respectively.
However, when it comes to nuclear spaces, things are getting much easier.
\begin{Tm}
Let $E$ be a locally convex Hausdorff space. Then, $E$ is nuclear if
and only if for every locally convex Hausdorff space $F$, $E \hat
\otimes_\pi F = E \hat \otimes_\epsilon F$.
\end{Tm}
A proof can be found in \cite[Theorem 50.1 p. 511]{Treves67}.
Thanks to this last theorem, we simply denote
\[
E \hat \otimes F \stackrel{\rm def.}{=} E \hat \otimes_\pi F = E
\hat \otimes_\epsilon F
\]
when one of the spaces $E$ or $F$ is nuclear. We also denote the
usual tensor product of two Hilbert spaces $E$ and $F$ by $E
\otimes F$. We recall the following result of Grothendieck on
tensor products, see \cite{Groth55}:
\begin{Pn}
Let $E$ be a complete locally convex space of functions defined
on a set $T$, such that its topology is finer than the pointwise
convergence topology, and assume $E$ to be nuclear. Then for every
complete locally convex space $F$, the tensor product $E
\hat{\otimes} F$ can be interpreted as the space of all functions
$f:T \to F$ such that for all $y^\prime \in F^\prime$, $t \mapsto
\langle y^\prime,f(t) \rangle _{F',F}$ is a function of $E$.
\end{Pn}

For $A,B \subseteq \mathbb N$, let $a:\ell_A \to \mathbb R$ and $b:\ell_B \to \mathbb R$
be two positive functions, such that the associated countably
Hilbert spaces $\mathcal F_{a}=\bigcap_{p\in \mathbb{N}} \ell^2_{a^p}$ and $\mathcal F_{b}
=\bigcap_{p\in \mathbb{N}} \ell^2_{b^p}$ are nuclear.
Since $\mathcal F_{a}$ is a complete locally convex space of functions
defined on the free commutative monoid $\ell_A$,
and since its topology is finer than the pointwise topology, that is since
$$|{\varphi_\alpha}_\lambda - {\varphi_\alpha}|^2  a_\alpha^{p}
 \leq \sum_{\alpha \in \ell_A} |{\varphi_\alpha}_\lambda - {\varphi_\alpha}| ^2
 a_\alpha^{p},
$$
$\mathcal F_{a} \hat{\otimes} \mathcal F_{b}$
 can be interpreted as the
 space of all elements of the form $\psi_{\beta,\alpha}$ such that for all
 $f=(f_\beta)_{\beta \in \ell_B} \in \mathcal F_{b}'$, $( \langle f,
 \psi_{\beta,\alpha} \rangle _{\mathcal F_{b}',
 \mathcal F_{b}})_{\alpha \in \ell_A} \in \mathcal F_{a}$.

\begin{Tm} \label{Tm:Tensor1}
It holds that
\[
\mathcal F_{a} \hat{\otimes} \mathcal F_{b} = \bigcap_{p,q}
\ell^2_{a^p} {\otimes} \ell^2_{b^q}.
\]
\end{Tm}
\begin{proof}
Let $\psi=( \psi_{\beta,\alpha} )$ be an element of $\mathcal F_{a} \hat{\otimes}
\mathcal F_{b}$ with $\psi_{\beta,\alpha}>0$ for any $\alpha \in \ell_A, \beta \in \ell_B$.
Then for all $f=(f_\beta)_{\beta \in \ell_B}
  \in \mathcal F_{b}'$, $( \langle f,\psi_{\beta,\alpha}
  \rangle _{\mathcal F_{b}',\mathcal F_{b}})_{\alpha \in \ell_A} \in
   \mathcal F_{a}$. In particular,
   we may choose $f_\beta=b_\beta^{\frac{q}{2}}$ for any $q \in \mathbb{N}$.
   Therefore,
\[
\begin{split}
\sum_{\alpha \in \ell_A} \sum_{\beta \in \ell_B} |\psi_{\beta,\alpha}|^2
 b_\beta^{q} a_\alpha^p
& \leq \sum_{\alpha \in \ell_A}\left| \sum_{\beta \in \ell_B} {b_
\beta^{\frac{q}{2}}}\psi_{\beta,\alpha}  \right|^2  a_\alpha^p \\
&=\sum_{\alpha \in \ell_A} \left| \langle (b_\beta^{\frac{q}{2}}),
(\psi_{\beta,\alpha}) \rangle _{\mathcal F_{b}',\mathcal F_{b}} \right|^2
a_\alpha^p < \infty
\end{split}
\]
Thus, $( \psi_{\beta,\alpha} ) \in  \bigcap_{p,q}  \ell^2_{a^p}
 {\otimes} \ell^2_{b^q}$.
In case $( \psi_{\beta,\alpha} ) \in \mathcal F_{a} \hat{\otimes} \mathcal F_{b}$
is an arbitrary element, then applying the last inequality to its positive real part, negative real
part, positive imaginary part and negative imaginary part yields the
requested result.\\

To prove the opposite direction, take $( \psi_{\beta,\alpha} ) \in  \bigcap_{p,q}  \ell^2_{a^p}
 {\otimes} \ell^2_{b^q}$. Then for all
 $f=(f_\beta)_{\beta \in \ell_B} \in \mathcal F_b'$
there exists $q$ such that $f \in \ell^2_{b^{-q}}$. Thus,
\[
\begin{split}
\left| \langle  f,(\psi_{\beta,\alpha}) \rangle _{\mathcal F_{b}',
\mathcal F_{b}} \right|^2
&=\left| \sum_{\beta \in \ell_B}
\overline{f_\beta}\psi_{\beta,\alpha}  \right|^2  \\
&=\left| \sum_{\beta \in \ell_B} \overline{f_\beta}
b_\beta^{-\frac q2}\psi_{\beta,\alpha}   b_\beta^{\frac q2}\right|^2  \\
&\leq  \left(\sum_{\beta \in \ell_B} |f_\beta|^2 b_\beta^{-q}\right)
\left(\sum_{\beta \in \ell_B} |\psi_{\beta,\alpha}|^2
 b_\beta^{q}\right)  \\
&= \|f\|_{\ell^2_{b^{-q}}}^2 \sum_{\beta \in \ell_B} |\psi_{\beta,\alpha}|^2
b_\beta^{q}
\end{split}
\]
Hence for all $p \in \mathbb{N}$,
\[
\sum_{\alpha \in \ell_A} \left| \langle  f,
(\psi_{\beta,\alpha}) \rangle _{\mathcal F_{b}',
\mathcal F_{b}} \right|^2  a_\alpha^p
\leq  \|f\|_{\ell^2_{b^{-q}}}^2 \sum_{\alpha \in \ell_A} \sum_{\beta \in \ell_B}
|\psi_{\beta,\alpha}|^2 b_\beta^{q}  a_\alpha^p  < \infty
\]
and so  $( \psi_{\beta,\alpha} ) \in \mathcal F_{a} \hat{\otimes} \mathcal F_{b}$.
\end{proof}

\begin{Pn}\label{Pn:Tensor2}
It holds that
\[
\bigcap_{p,q} \ell^2_{a^p} {\otimes} \ell^2_{b^q} =
 \bigcap_{p} \ell^2_{a^p} {\otimes} \ell^2_{b^p}
\]
\end{Pn}
\begin{proof}
One direction is clear. The other direction follows from the inclusion
\[
\ell^2_{a^{\max\{p,q\}}} {\otimes} \ell^2_{b^{\max\{p,q\}}}
 \subseteq  \ell^2_{a^p} {\otimes} \ell^2_{b^q}
\]
\end{proof}

Applying Theorem \ref{Tm:Tensor1} and Proposition \ref{Pn:Tensor2} ,we obtain:
\begin{equation}
\label{Eq:Tensor3}
\mathcal F_{a} \hat{\otimes} \mathcal F_{b}=
\bigcap_{p \in \mathbb{N}}
\left\{ (\psi_{\alpha,\beta}) :
 \sum_{(\alpha,\beta) \in \ell_A \times \ell_B}
|\psi_{\beta,\alpha}|^2 (b_\beta a_\alpha)^p <\infty \right\}.
\end{equation}
Now, we can concatenate the indices in an obvious way. We define
$C=2B \cup (2A-1)$ (disjoint union), $P_B:C \to B$, $P_A:C \to A$
the appropriate projections,
$c_\gamma=b_{P_B(\gamma)}a_{P_A(\gamma)}$,
and
$\psi_\gamma=\psi_{{P_B(\gamma)},{P_A(\gamma)}}$. Therefore, we
may write
\[
\mathcal F_{a} \hat{\otimes} \mathcal F_{b}=\bigcap_{p \in \mathbb{N}}
\left\{ ( \psi_{\gamma})_{\gamma \in \ell_C} : \sum_{\gamma \in \ell_C}
|\psi_{\gamma}|^2  c_\gamma^p <\infty \right\},
\]
and
\[
(\mathcal F_{a} \hat{\otimes} \mathcal F_{b})'=\bigcup_{p \in \mathbb{N}}
\left\{ ( \psi_{\gamma})_{\gamma \in \ell_C} : \sum_{\gamma \in \ell_C}
|\psi_{\gamma}|^2  c_\gamma^{-p} <\infty \right\}.
\]
\begin{Pn}\label{Pn:Tensor4}
Let $a:\ell_A \to \mathbb R$ and $B:\ell_B \to \mathbb R$ be two admissible positive functions.
Then:
\begin{enumerate}[(a)]
\item $c:C \to \mathbb R$ (where $C=2B\cup(2A-1)$ and $c_\gamma=c_{P_B(\gamma)}a_{P_A(\gamma)}$)
 is admissible.
\item If $a$ and $b$ are both $d$-regular, then so is $c$.
\item If $a$ and $b$ are both superexponential, then so is $c$.
\end{enumerate}
\end{Pn}
\begin{proof}[Proof]
$c$ is admissible, since $c_0 =b_0a_0=1$
and $c_{e_n}>1$ for all $n \in C$. Moreover,
\[
\sum _{n \in C} \frac 1{c_{e_n}^d-1}=\sum _{n \in A} \frac 1{a_{e_n}^d-1}+
\sum _{n \in B} \frac 1{b_{e_n}^d-1},
\]
and hence $d$-regularity of both $a$ and
$b$ yields $d$-regularity of $c$.
Finally, if both $a$ and $b$ are superexponential, clearly so is $c$.
\end{proof}
Finally, we give the following theorem.
\begin{Tm}\label{Tm:DualTensor}
Let $E$, $F$ be two Fr\'{e}chet spaces. If $E$ is nuclear,
we have the canonical isomorphism
\[
(E \hat{\otimes} F)'=E' \hat{\otimes} F'
\]
\end{Tm}
A proof, in case both $E$ and $F$ are nuclear, is given in
\cite[(50.19), p. 525]{Treves67}. We can now state:
\begin{Tm}
A tensor product of two V{\aa}ge spaces is also a V{\aa}ge space.
\end{Tm}
\begin{proof}[Proof]
Applying \eqref{Eq:Tensor3}, Proposition \ref{Pn:Tensor4} and Theorem \ref{Tm:VageSpace},
$(\mathcal F_{a} \hat{\otimes} \mathcal F_{b})^\prime$ is a V\aa ge space. Theorem \ref{Tm:DualTensor} yields the requested result.
\end{proof}


\section{An extension of the space of tempered distributions}
\label{sec:schwartz}
\setcounter{equation}{0}
In this section we consider the special case
$\ell_A=\mathbb{N}_0$ (i.e. $A=\{1\}$), and
\[
a_n=(n+1)^2. 
\]
The corresponding space $\mathcal F_{a}$ (defined by
\eqref{eq:Fdef}) is identified with the Schwartz space $ \mathscr
S$ of rapidly decreasing smooth functions, and its dual is the
space $\mathscr S^\prime$ of tempered distributions. We will show
(see Proposition \ref{pn:Schwartz} below) that $\mathscr
S^\prime$ is a regular admissible space, but it is not a V\aa ge
space. We will also construct a V\aa ge space containing
$\mathscr S^\prime$.\\

We recall (see \cite[Chapter IV, Section 2, p. 303]{sansone},
\cite[p. 105]{MR0372517})
 that the
{\sl Hermite polynomials} $h_n(x)$ are
defined by
\begin{equation}
\label{hermitepol}
h_n(x) =
(-1)^ne^{x^2}\frac{d^n}{dx^n}\left(e^{-x^2}\right),\quad n=0,1,\ldots,
\end{equation}
Various notations and conventions are given for these polynomials. See
the discussion on the end of page 105 of \cite{MR0372517}. In particular,
the multiplicative factor $(-1)^n$ (which does not appear in Sansone's book
\cite{sansone}) insures that the factor of $x^n$ in $h_n$ is positive.
The {\sl Hermite functions} $\xi_n(x)$ are defined by
\begin{equation}
\label{hermitefunction}
\xi_n(x) = \pi^{-\frac14}(2^nn!)^{-\frac12}e^{-\frac12 x^2}
h_{n}(x),\quad n=0,1,2,\ldots.
\end{equation}
The Hermite functions $(\xi_n)_{n \in \mathbb{N}_0}$ form an
orthonormal basis of $\mathbf L_2(\mathbb{R},dx)
$. The {\sl Schwartz space}
$\mathscr S$ of smooth rapidly decreasing functions on
$\mathbb{R}$ is defined by
\[
\mathscr S=\left\{ f \in C^{\infty}(\mathbb{R}) : \sup_{x \in
\mathbb{R}} \left|x^pf^{(q)}(x)\right|<\infty \text { for all }
p,q \in \mathbb{N}_0 \right\}
\]
The Hermite functions are elements in the Schwartz space, and we
have (see \cite[Theorem V.13 p. 143]{MR751959}):
\[
\mathscr S=\left\{ f=\sum_{n \in \mathbb{N}_0}f_n \xi_n : \sum_{n
\in \mathbb{N}_0}|f_n|^2(n+1)^{2p}<\infty  \text{ for all } p \in
\mathbb{N}\right\}.
\]
Identifying $\sum_{n \in \mathbb N_0}f_n \xi_n$ with $(f_n)_{n \in \mathbb N_0}$
allows to identify
$\mathbf L_2(\mathbb{R},dx)$ with $\ell^2(\mathbb N_0)$ and $\mathscr S$ with
$\mathcal F_{a}$.

\begin{Pn}
\label{pn:Schwartz}
$\mathscr S^\prime$ is a regular admissible space which is nuclear, but is not
a V{\aa}ge space.
\end{Pn}

\begin{proof}[\indent Proof]
First, we note that defining $a_n=(n+1)^2$ implies that
$a:n \mapsto a_n$ is a 1-regular admissible function. It
is indeed admissible since $a_0=1$ and $a_1=4>1$, and it is
indeed 1-regular, since the sum in \eqref{eq:ineq1} is over the
finite set $A=\{1\}$ and in particular does converge.
Therefore, $\mathscr
S^\prime$ is a 1-regular admissible space. Since, $\sum _{n \in
\mathbb{N}_0} \left((n+1)^2\right)^{-1} < \infty$, $\mathscr S$
is nuclear, and hence $\mathscr S^\prime$ is also nuclear (of
course, the nuclearity of $\mathscr S$ and $\mathscr S^\prime$ is
a standard result; see for instance \cite{Treves67}). Since
$a$ is not superexponential, that is
\[
(n+1)^2(m+1)^2 \not \leq (n+m+1)^2,
\]
$\mathscr S^\prime$ is not a V\aa ge space.
\end{proof}

We define the following subspace of $\mathscr S$
\[
\mathscr G = \left\{\sum_{n=0}^\infty f_n \xi_n :   \sum_{n
\in \mathbb{N}_0}|f_n|^2 2^{np}<\infty  \text{ for all } p \in
\mathbb{N}\right\}.
\]
\begin{Pn}
$\mathscr G^\prime$ is a V\aa ge space containing the Schwartz
space $\mathscr S^\prime$ of tempered distributions.
\end{Pn}
\begin{proof}[\indent Proof]
Using the identification of $\sum_{n \in \mathbb N_0}f_n \xi_n$
with $(f_n)_{n \in \mathbb N_0}$, and defining $a_n=2^n$, we have that $\mathscr G$ is the corresponding countably
Hilbert space $\mathcal F_{a}$ associated $a$ (and as before,
$\mathbf L_2(\mathbb{R},dx)$ is identified with $\ell^2(\mathbb N_0)$).
Clearly, $a$ is a 1-regular admissible
function. It is indeed admissible since $a_0=1$ and $a_1=2>1$,
and it is indeed 1-regular, since the sum in \eqref{eq:ineq1} is
over the finite set $A=\{1\}$ and in particular does converge.
Therefore, $\mathscr G^\prime$ is a 1-regular
admissible space. Since $a:n \mapsto 2^n$ is an
exponential function, $\mathscr G^\prime$ is a V\aa ge space and
is in particular nuclear. Moreover, the natural embeddings
$\mathscr G \subseteq \mathscr S$ and $\mathscr S' \subseteq
\mathscr G'$ are clearly continuous. Hence $\mathscr G$ is a
closed subspace of $\mathscr S$, and $\mathscr S'$ is a closed
subspace of $\mathscr G'$.
\end{proof}

\begin{Tm}
\label{thm:entire}
$\mathscr G$ is the space of all entire functions $f(z)$ such that
\[
 \iint_\mathbb{C} \left|f(z)\right|^2 e^{\frac{1-2^{-p}}{1+2^{-p}}x^2
 -\frac{1+2^{-p}}{1-2^{-p}}y^2}dxdy<\infty \quad \text { for all } p \in \mathbb N.
\]
\end{Tm}

In the proof we make use of two results of Hille. The first result
appear in  \cite[formula (1.3), p. 81]{MR0000871} and
\cite[Theorem 2.2 p. 885]{Hille_Duke}. For the second formula, see
\cite[formula (2.1) p. 82]{MR0000871} and \cite[p.
439-440]{MR1502747}. In that last paper one can also find a
history of the formula.

\begin{Tm} (Hille, \cite{MR0000871})
\label{Tm:Hille1} The domain of absolute convergence of the
series $\sum_{n=0}^\infty F_n \xi_n(z)$ is the strip
$S_\tau=\left\{z \in \mathbb{C} : |{\rm Im}(z)| <\tau \right\}$,
where
\[
\tau=-\limsup_{n \to \infty}(2n+1)^{-\frac12}\log|F_n|.
\]
\end{Tm}

\begin{Tm} (Hille, \cite{MR1502747})
\label{Tm:Hille2}
The series  $\sum_{n=0}^\infty{\xi_n(u)\xi_n(v)s^n}$ converges for
arbitrary complex values of $u$ and $v$ when $|s|<1$, and
\[
\sum_{n=0}^\infty{\xi_n(u)\xi_n(v)s^n}=\pi^{-\frac 12}(1-s^2)^{-\frac 12}
e^{-\frac{(1+s^2)(u^2+v^2)-4svu}{2(1-s^2)}}.
\]
\end{Tm}

Furthermore, we make use of the easy following proposition. See
\cite[\S 6, p. 61]{MR1478165}. The space $\Gamma_\alpha$ bears
various names, and in particular is called the Fock space.
\begin{Pn}
\label{Pn:Alpay}
For all $0<\alpha \leq 1$,
\[
\Gamma_\alpha=\left\{f \text{ is entire  }:
 \frac{\alpha}{\pi}\iint_\mathbb{C}\left|f(z) \right|^2
 e^{-\alpha|z|^2}dxdy<\infty   \right\}
\]
is a Hilbert space with a reproducing kernel
$K_\alpha(z,w)=e^{\alpha \overline{w}z}$.
\end{Pn}

\begin{proof}[\indent Proof
 of Theorem \ref{thm:entire}]
Let $p \in \mathbb{N}$. For each $f \in
\mathscr G$, $f= \sum_{n=0}^\infty f_n \xi_n$ whereas $\sum_{n=0}^\infty |f_n|^2
2^{np}<\infty$ for all $p \in \mathbb{N}$. In particular, $\lim_{n \to \infty}|f_n|^2
 2^{n}=0$.
Therefore, for every $n$ large enough, $\log |f_n|<-n\log\sqrt2$.
Thus, in the notations of Theorem \ref{Tm:Hille1},
\[
\tau=-\limsup_{n \to \infty}(2n+1)^{-\frac12}\log|f_n| \geq \liminf_{n
\to \infty}(2n+1)^{-\frac12}n\log\sqrt2=\infty.
\]
Therefore, denoting by $\mathscr G_p$ the space of all functions
 $f=\sum_{n=0}^\infty f_n\xi_n$ subject to $\sum_{n=0}^\infty |f_n|^22^{np}<\infty$ (i.e., $\mathscr G_p \cong \ell^2_{a^p}$),
 it is a Hilbert space of entire functions,
and in particular, $\mathscr G=\bigcap_{p \in \mathbb{N}}
\mathscr G_p$ is a countably Hilbert space of entire functions.\\

Now, since $\left( \xi_n2^{-\frac{np}{2}} \right)_{n \in \mathbb{N}_0}$ is
an orthonormal basis for $\mathscr G_p$,
and denoting $s=2^{-p}$,
the reproducing kernel of the Hilbert space $\mathscr G_p$ is given by
\[
G(z,w)=\sum_{n=0}^\infty \xi_n(z) \overline{\xi_n(w)} 2^{-np}=
\sum_{n=0}^\infty \xi_n(z) \xi_n(\overline{w}) s^{n}.
\]
Applying Theorem \ref{Tm:Hille2}, we have that
\[
G(z,w)=\pi^{-\frac 12}(1-s^2)^{-\frac 12}e^{-\frac{(1+s^2)
(z^2+\overline{w}^2)-4sz\overline{w}}{2(1-s^2)}}
\]
Denoting $r(z)=\pi^{-\frac 14}(1-s^2)^{-\frac 14}e^{-\frac {1+s^2}{2(1-s^2)}z^2}$,
considering the kernel $K_\alpha(z,w)=e^{\alpha \overline{w}z}$
with its associated Hilbert space $\Gamma_\alpha$ for
$\alpha=\frac{2s}{1-s^2}$ (see Proposition \ref{Pn:Alpay}),
we have that
\[
G(z,w)=r(z)K_\alpha(z,w)r(\overline{w}).
\]
Therefore, the space $\mathscr G_p$ is equal to the space of functions of the form $f=rg$, with
$g\in \Gamma_\alpha$ and norm
\[
\|f\|_{\mathscr G_p}=\|g\|_{\Gamma_\alpha}.
\]
We note that for $z=x+iy$,
\[
\begin{split}
2 \cdot \frac{1+s^2}{2(1-s^2)}{\rm Re}(z^2) - \frac{2s}{1-s^2}|z|^2
&=\frac{(1+s^2)(x^2-y^2)}{1-s^2}-\frac{2s(x^2+y^2)}{1-s^2}\\
&=\frac{1-s}{1+s}x^2-\frac{1+s}{1-s}y^2.
\end{split}
\]
Thus,
\[
|r(z)|^{-2}=\sqrt{\pi(1-s^2)}e^{\frac{1-s}{1+s}x^2-\frac{1+s}{1-s}y^2},
\]
and, with $K_p=\frac{2^{1-p}}{\sqrt{\pi(1-2^{-2p})}}$,
\[
\mathscr G_p=\left\{ f \text{ is entire  }:\, \|f\|_{\mathscr G_p}^2=K_p
 \iint_\mathbb{C} \left|f(z)\right|^2 e^{\frac{1-2^{-p}}{1+
 2^{-p}}x^2-\frac{1+2^{-p}}{1-2^{-p}}y^2}dxdy<\infty  \right\},
\]
and in particular $\mathscr G=\bigcap_{p \in \mathbb N}\mathscr G_p$ is
the space of all entire functions $f(z)$ subject to
\[
 \iint_\mathbb{C} \left|f(z)\right|^2 e^{\frac{1-2^{-p}}{1+2^{-p}}x^2
 -\frac{1+2^{-p}}{1-2^{-p}}y^2}dxdy<\infty \quad \text { for all } p \in \mathbb N.
\]
\end{proof}

We note that there are functions in the Schwartz space $\mathscr
S$  which have no analytic continuation to an entire function but only
to a holomorphic function on some strip. For
example, one may consider the function
\[
F(x)=\sum_{n \in \mathbb{N}_0} e^{-\sqrt{2n+1}}\xi_n(x) \in \mathscr S
\]
It is indeed in the Schwartz space since $\sum_{n \in
 \mathbb{N}_0}e^{-2\sqrt{2n+1}}(n+1)^{2p}<\infty$ for all $p
 \in \mathbb{N}$. However,
\[
\tau=-\limsup_{n \to \infty}(2n+1)^{-\frac12}\log|F_n| =
\liminf_{n \to \infty}(2n+1)^{-\frac12}(2n+1)^{\frac 12}=1.
\]
Furthermore, clearly there are functions in the Schwartz space
which have no analytic continuation to any holomorphic function, e.g.,
functions with compact support.

\begin{Rk}
As mentioned in the introduction,
connections with the theory of hyperfunctions will be considered elsewhere.
\end{Rk}

\begin{Rk}
Having now $\mathscr G^\prime$ at hand, one can
consider the tensor product $\mathscr G^\prime \otimes
\bigS_{-1}$, where $\bigS_{-1}$ is the Kondratiev space of stochastic
distributions (see Section \ref{sec:7}). By Section \ref{sec:5} it is a V\aa ge space,
and can be an appropriate setting to study stochastic
linear systems. This will be done in a future publication.
\label{rk:rk}
\end{Rk}

\begin{Rk}
If we define  $a_0=1$ and for $n>0$ $a_n=2^{2^n}$, then $a$
is again a $1$-regular admissible function. The associated
countably Hilbert space
$\mathcal F_{a}$ is a space of entire functions.
Moreover, since the function $a$ is superexponential
its dual is clearly a V\aa ge space.
We can define many other V\aa ge spaces in a similar manner.
\end{Rk}

\section{The Kondratiev spaces}
\setcounter{equation}{0}
\label{sec:7}
In this section we consider the special case $\ell_A=\ell$ (i.e.
$A=\mathbb N$), and
\[
a_\alpha=(2\mathbb N)^\alpha.
\]
Then the
corresponding space $\mathcal F_{a}$ (defined by
\eqref{eq:Fdef}) is identified with the Kondratiev space of
Gaussian test functions $\bigS_1$, and its dual is the Kondratiev
space of Gaussian stochastic distributions $\bigS_{-1}$. We will
show (see Proposition \ref{pn:Schwartz} below) that
$\bigS_{-1}$ is a V\aa ge space. We also consider the Kondratiev space of
Poissonian stochastic distributions, and show that it is also a V\aa ge space.\\

We first need to recall a few definitions pertaining to the white noise space.
The function $s \mapsto e^{-\frac12 \|s\|^2_{\mathbf L_2(\mathbb{R},dx)}}$
is positive definite on the
Schwartz space of real-valued functions $\mathscr S_\mathbb R$, and continuous at the origin.
The Bochner-Minlos theorem (see for instance \cite[p. 10-11]{MR2105995})
insures the existence of a probability measure $d\mu$ on the Borel $\sigma$-algebra of
the dual space $\mathscr S_\mathbb R^\prime$, such that
\[
e^{-\frac12 \|s\|^2_{\mathbf L_2(\mathbb{R},dx)}}
=\int_{\mathscr S_\mathbb R^\prime}e^{i\langle s', s \rangle}d \mu(s')
\quad \text{ for all } s \in \mathscr S_{\mathbb R},
\]
where the brackets denote the duality between $\mathscr S_\mathbb R$ and
$\mathscr S_\mathbb R^\prime$. This equality induces an isometric map
\[
s\mapsto Q_s,\quad{\rm where}\quad Q_s(s')=\langle s',s\rangle \quad (s \in \mathscr S_{\mathbb R}, s' \in \mathscr S'_{\mathbb R})
\]
from $\mathscr S_\mathbb R\subset \mathbf L_2(\mathbb R,dx)$ into
$\mathbf L_2(\mathscr S_\mathbb R^\prime,\mathcal B, \mu)$.\\

The space $\mathcal W= \mathbf L_2(\mathscr S_\mathbb R^\prime,\mathcal B, \mu)$
is called the Gaussian white noise space.
We recall that the
{\sl Hermite polynomial functionals} $(H_\alpha)_{\alpha \in \ell}\subseteq \mathcal W$, which are defined by
\[
H_\alpha(s')=\prod_{k=1}^\infty h_{\alpha_k}(Q_{{\xi_{k-1}}}(s'
)) \quad (s' \in \mathscr S_\mathbb R^\prime),
\]
form an orthogonal basis of $\mathcal W$,
where $(h_k)$ and $(\xi_k)$ denote respectively the
Hermite polynomials and the Hermite functions (see \eqref{hermitepol}
and \eqref{hermitefunction}).
 More precisely,
\begin{equation}
\mathcal W=\left\{ \sum_{\alpha \in \ell} f_\alpha H_\alpha: \sum_{\alpha \in \ell}
|f_{\alpha}|^2\alpha!< \infty \right\},
\label{eq:gauss}
\end{equation}
where $\alpha!=\prod_{k=1}^\infty (\alpha_k)!$.
The Kondratiev space of Gaussian test function $\bigS_1$ is defined by
\[
\bigS_1=\left\{ \sum_{\alpha \in \ell} f_\alpha H_\alpha: \sum_{\alpha \in \ell}
 |f_{\alpha}|^2(2\mathbb N)^{\alpha p}(\alpha!)^2< \infty \text { for all } p
 \in \mathbb N \right\}.
\]
Clearly, $\mathcal W$ can be identified with $\ell^2$ using the isometry
\[
\sum_{\alpha \in \ell} f_\alpha H_\alpha \mapsto \left(f_\alpha(\alpha!)^{\frac12}
\right)_{\alpha\in\ell_A}.
\]
In a similar way,  defining
\[
\bigS_{1,p}=\left\{ \sum_{\alpha \in \ell} f_\alpha H_\alpha: \sum_{\alpha \in \ell}
 |f_{\alpha}|^2(2\mathbb N)^{\alpha p}(\alpha!)^2< \infty \right\},
\]
it can be identified with $\ell^2(a^p)$ (for $a_\alpha=(2\mathbb N)^\alpha$).
Hence,
\[
\mathcal F_a = \bigcap_{p\in \mathbb N} \ell^2(a^p) \cong \bigcap_{p\in \mathbb N} \bigS_{1,p} =\bigS_{1}.
\]

The {\sl Wick product} of two formal series
$f=\sum_{\alpha \in \ell_A}f_\alpha H_\alpha$ and $g=\sum_{\alpha \in \ell_A}g_\alpha H_\alpha$,
denoted by $\diamond$,  was introduced by Hida and Ikeda, see \cite{Hida_Ikeda}.
It is defined by
\begin{equation}
f \diamond g =\sum_{\gamma \in \ell} \left( \sum_{\alpha+\beta=\gamma}{f_\alpha g_\beta} \right) H_\gamma.
\end{equation}
The Kondratiev space of Gaussian stochastic distributions $\bigS_{-1}$, which is the dual of $\bigS_1$ and can be defined by
\[
\bigS_{-1}=\left\{ \sum_{\alpha \in \ell} f_\alpha H_\alpha: \sum_{\alpha \in \ell}
 |f_{\alpha}|^2(2\mathbb N)^{-\alpha p}< \infty \text { for some } p
 \in \mathbb N \right\} \cong \mathcal F_a',
\]
is not only closed under the Wick product, but as we present in the following proposition, also a V\aa ge space.
The following proposition is a result of V\aa ge, proved in 1996, see
\cite{vage96b} and \cite[p. 118]{MR1408433}.
\begin{Pn}
The Kondratiev space of Gaussian stochastic distributions $\bigS_{-1}=\bigS_1'$ is a V\aa ge space.
\end{Pn}
\begin{proof} [\indent Proof]
We note that $a: \alpha \mapsto (2\mathbb N)^\alpha$ is a $2$-regular admissible positive function.
It is admissible since $a_0=(2n)^0=1$ and
$a_{e_n}=2n>1$, and it is $2$-regular since $\sum _{n \in \mathbb{N}}
\frac 1{(2n)^2-1}<
\infty$. Therefore,$\bigS_{-1} \cong \mathcal F_{a}'$  is a $2$-regular admissible
space. Since $a$ is exponential, $\bigS_{-1}$ is a
V\aa ge space.
\end{proof}

Hida's theory can also be applied to Poisson processes. One then considers
the function
\[
\exp \left[\int_{\mathbb R}(e^{i s(x)}-1)dx\right] \quad(s\in\mathscr S_{\mathbb R})
\]
which is positive definite on $\mathscr S_\mathbb R$, and continuous at the
origin.
Here too, the Bochner-Minlos theorem insures the existence of a probability
measure $\pi$ on $\mathscr S^\prime_\mathbb R$ and such that
\[
\exp \left[\int_{\mathbb R}(e^{i s(x)}-1)dx\right]
=\int_{\mathscr S_\mathbb R^\prime}e^{i\langle s', s \rangle}d \pi(s').
\]

The Poissonian white noise space is $\mathcal W^{\pi}= \mathbf L_2
(\mathscr S_\mathbb R^\prime,\mathcal B, \pi)$, and admits a representation
of the form \eqref{eq:gauss},
replacing the Hermite polynomial functionals
$(H_\alpha)_{\alpha \in \ell}$ with the Charlier polynomial functionals
$(C_\alpha)_{\alpha \in \ell}$, computed in terms of the Poisson-Charlier
polynomials; see \cite[p. 185]{MR1408433}; we refer to \cite[Chapter II,
\S 2.81, p. 34-35]{MR0372517} for the Poisson-Charlier polynomials.
More precisely,
\[
\mathcal W^{\pi}=\left\{ \sum_{\alpha \in \ell}
f_\alpha C_\alpha: \sum_{\alpha \in \ell} |f_{\alpha}|^2\alpha!< \infty \right\}.
\]
The Kondratiev space of Poissonian test function $\bigS_1^\pi$ is defined by
\[
\bigS_1^\pi=\left\{ \sum_{\alpha \in \ell} f_\alpha C_\alpha:
\sum_{\alpha \in \ell} |f_{\alpha}|^2(2\mathbb N)^{\alpha p}(\alpha!)^2< \infty
\text { for all } p \in \mathbb N \right\}.
\]
Since the associated positive function $a:\alpha \mapsto a_\alpha$ remains the same as
before, we conclude the following proposition.
\begin{Pn}
The Kondratiev space of Poissonian stochastic distributions is a V\aa ge space.
\end{Pn}

We refer to \cite[Theorem 3.7 p. 192]{MR1851117} for another example of a space where V\aa ge inequality holds in the setting of white noise space analysis.

\section{State space theory and V\aa ge spaces}
\setcounter{equation}{0} \label{sec:8}
The results presented in \cite{alp,aa_goh} for the case of the
Kondratiev space $\mathcal S_{-1}$ of stochastic distributions
extend to general V\aa ge spaces. We refer to  \cite{MR2584037,Fuhrmann,MR0255260}
for general background on the theory
of linear systems when the coefficient space is $\mathbb C$ (we also refer to
\cite{MR2002b:47144} for a survey),
and to the papers \cite{SSR,MR1071708}
and to the book \cite{MR839186}
for more information on linear system on commutative rings, and in particular
for the notions of controllable and observable pairs. These various notions
are also reviewed in \cite{alp}.\\

We begin this section with the following proposition. For
completeness we give a short outline of the proof.
\begin{Pn}
A matrix-valued rational function $R(z)=p(z)(q(z))^{-1}$ for which
$E(q(0))\not=0$ can be written as
\begin{equation}
\label{eq:rat}
R(z)={D}+z{C}(I-z{A})^{-1}{B}
\end{equation}
where $A$, $B$, $C$, and $D$ are
matrices of appropriate dimensions and with entries in the ring $\bigR$.
\label{pn:rational}
\end{Pn}

{\bf Proof:} We first note that a constant and the function $z$
trivially have realizations of the form \eqref{eq:rat}.
Furthermore, if $R_1$ and $R_2$ are in the form \eqref{eq:rat},
\[
R_j(z)=D_j+C_j(I_{N_j}-zA_j)^{-1}B_j,\quad j=1,2,
\]
one has the realization formulas
\begin{equation}
R_1(z)R_2(z)=D+C(zI_N-A)^{-1}B, \label{real_mul}
\end{equation}
where $N=N_1+N_2$, $D=D_1D_2$ and
\[
C=\begin{pmatrix}C_1& D_1C_2\end{pmatrix},\quad
B=\begin{pmatrix}B_1D_2\\ B_2\end{pmatrix}\quad{\rm and}\quad
A=\begin{pmatrix}A_1 &B_1C_2\\0&A_2\end{pmatrix},
\]
and
\begin{equation}
\label{real_sum}
R_1(z)+R_2(z)=D+C(I_N-zA)^{-1}B,
\end{equation}
where $N=N_1+N_2$, $D=D_1+D_2$ and
\[
C=\begin{pmatrix}C_1& C_2\end{pmatrix},\quad
B=\begin{pmatrix}B_1\\ B_2\end{pmatrix}\quad{\rm and}\quad
A=\begin{pmatrix}A_1 &0\\0&A_2\end{pmatrix}.
\]
for the product and sum (provided the dimensions of $R_1$ and
$R_2$ are such that these make sense). Next, if $R$ is $\mathcal
R^{p\times p}$ valued of the form \eqref{eq:rat} and $D$ is
invertible in $\mathcal R^{p\times p}$ then,
\[
R^{-1}(z)=D^{-1}-D^{-1}C(I_N-z(A-BD^{-1}C)^{-1}BD^{-1},
\]
which is also of the form \eqref{eq:rat}. To conclude it remains
to verify that realization is a property which stays under
concatenation: If $R_1$ and $R_2$ are of the form \eqref{eq:rat}
so are the functions
\[
\begin{pmatrix}R_1  \\ R_2\end{pmatrix}\quad{\rm and}\quad
\begin{pmatrix}R_1  & R_2\end{pmatrix},
\]
provided the dimensions make sense.\mbox{}\qed\mbox{}\\

We note that, as was already mentioned after Definition \ref{def4}, one can compute the value of a rational function of
the form \eqref{eq:rat1} at every point $f\in\bigR$ such that
$E(q(f))\not =0$. Let $q(z)=\sum_{m=0}^M  q_m z^{m}$. Then, this
last condition can be rewritten as
\[
\sum_{m=0}^M E(q_m)(E(f))^m\not=0.
\]
Similarly, one can compute
\eqref{eq:rat} at every  $f\in\mathcal R$ such that
$(I-E(fA))$ is invertible.\\

It is convenient to introduce the operators
$D_n$, $n=1,2,\ldots$ defined by
\[
D_n(x^\alpha)=
  \begin{cases}
   \alpha_n x^{\alpha - e_n} & \text{if } \alpha_n>0\,\\
   0 & \text{otherwise}
  \end{cases}
\]
and by linear extension to any finite linear combination of
such elements. We have in particular $D_n(XY)=D_n(X)Y+XD_n(Y)$.\\

As in \cite{alp} for $\mathcal S_{-1}$, given a V\aa ge space
$\mathcal R$ we define a rational function to be an expression of
the form
\[
D+zC(I-zA)^{-1}B,
\]
where $A,B,C$ and $D$ are matrices of
appropriate dimensions and with entries in $\mathcal R$. See Proposition
\ref{pn:rational} above.
Before giving a sample result we recall that a pair $( C,
A)\in{\mathcal R}^{p\times N}\times{\mathcal R}^{N\times N}$ is
called {\it observable} if the map
\[
 f\mapsto\begin{pmatrix}C f& CA f& C A^2
f&\cdots\end{pmatrix}
\]
is injective from ${\mathcal R}^N$ into $\left({\mathcal
R}^p\right)^{\mathbb N}$. See \cite[\S 2.2 p. 58]{MR839186}. In
\cite{alp} it is proved for $\mathcal S_{-1}$ that an equivalent condition is:

\[
C(I-z A)^{-1} f\equiv0_{{\mathcal R}}^{p\times N}\Longrightarrow
 f=0_{{\mathcal R}}^{N}.
\]
The proof is the same for any V\aa ge space.

\begin{Tm}
Let $\hat{h}$ be a rational function  with realization
\begin{equation}
\label{real:123} \hat{h}(z)=D+zC(I-zA)^{-1}B.
\end{equation}
If the realization $E[\hat{h}](z)=E[D]+zE[C](I-zE[A])^{-1}E[B]$
is observable, then the realization \eqref{real:123} is observable.
\label{tm:obs}
\end{Tm}

\begin{proof}[\indent Proof]
We assume that $A\in\mathcal R^{N\times N}$. Let $f \in \bigR^N$  be such that
$C(I-zA)^{-1}f\equiv 0$. We want to show that $f_\alpha=0$ for
all $\alpha \in \ell_A$. Since
\[
E[\hat{h}](z)=E[D]+zE[C](I-zE[A])^{-1}E[B]
\]
is an observable realization, we
have that $E[C](I-zE[A])^{-1}E[f]\equiv 0$ implies $E[f]=0$, and thus
$f_0=0$. Now, since
\[
\begin{split}
D_n(C(I-zA)^{-1}f)&=
D_n(C)(I-zA)^{-1}f+CD_n((I-zA)^{-1})f+\\
&\hspace{5mm}+
C(I-zA)^{-1}D_n(f)\\
&=0,
\end{split}
\]
and since $E[f]=0$, we have that $E[C](I-zE[A])^{-1}E[D_n(f)]=0$,
and thus $f_{e_n}=0$. Furthermore, by a simple induction, since
$$
D_n^m(C(I-zA)^{-1}f)=\sum_{k<m}{U_k D_n^k(f)}+C(I-zA)^{-1}D_n^m(f)
$$
for some $U_k$, and since $D_n^m(C(I-zA)^{-1}f)=0$ we have that
\[
E[C](I-zE[A])^{-1}E[D_n^m(f)]=0,
\]
and therefore $f_{m e_n}=0$.
Thus, $f_\alpha=0$ for all $\alpha \in \ell_A$ such that $\alpha=(0,
\dots, 0, \alpha_n, 0\dots)$.\\
We may complete this proof as in \cite{alp}.
\end{proof}

In conclusion, Theorem \ref{tm:obs} as well as Problem \ref{pb:NP} (or more precisely, its
solution presented in \cite{aa_goh}) suggest that most of
the classical linear system
theory can be extended to our setting. This is important when one wants to take into account
stochastic aspects of the theory. One such avenue consists of continuing the
line of research initiated in \cite{aa_goh,al_acap, alp}. One then considers input-output systems of the form
\[
y_n=\sum_{k=0}^n h_k u_{n-k},\quad n=0,1,\ldots
\]
where the input sequence $(u_n)_{n\in\mathbb N_0}$ and the impulse response
$(h_n)_{n\in\mathbb N_0}$ are in some V\aa ge space. The choice
of the given V\aa ge space
is done to express for instance that the system is stochastic
(then one
chooses $\mathcal S_{-1}$). When the sequences consist of complex numbers,
the product reduces to the product of complex numbers, and we are back in the
classical theory.\\


As was already mentioned in Remark \ref{rk:rk}, another avenue is to
define a stochastic linear system as a continuous mapping from
the nuclear space $\mathscr G\otimes \mathcal S_{1}$ into its
dual, to use Schwartz's kernel theorem and then to follow
Zemanian's approach to linear systems. Test functions are now
functions of the form (we write $\w$ rather than $s^\prime$ for
the variable in $\mathcal S_1$)
\begin{equation}
\label{test}
s(t,\w)=\sum_{\substack{n\in{\mathbb N}_0\\
\alpha\in\ell}}H_\alpha(\w)\xi_n(t)c_{n,\alpha}
\end{equation}
where the coefficients $c_{n,\alpha}$ are in ${\mathbb C}$ and
subject to
\[
\|s\|_{p,q}^2=\sum_{n=0}^\infty\sum_{\alpha\in\ell}
(\alpha!)^2|c_{n,\alpha}|^2(2\mathbb N)^{q\alpha}2^{np}<\infty,\quad \forall
p,q\in{\mathbb N}.
\]

See \cite{zemanian} for the latter. Since the dual of $\mathscr
G\otimes \mathcal S_{1}$ is a V\aa ge space, one can get more
precise results than
the ones in \cite{zemanian}.\\


\bibliographystyle{plain}

\begin{thebibliography}{10}

\bibitem{MR2002b:47144}
D.~Alpay.
\newblock {\em The {S}chur algorithm, reproducing kernel spaces and system
  theory}.
\newblock American Mathematical Society, Providence, RI, 2001.
\newblock Translated from the 1998 French original by Stephen S. Wilson,
  Panoramas et Synth\`eses. [Panoramas and Syntheses].

\bibitem{aa_goh}
D.~Alpay and H.~Attia.
\newblock An interpolation problem for functions with values in a commutative
  ring.
\newblock Operator Theory: Advances and Applications, vol. 218 (2011), p. 1-17.


\bibitem{acs1}
D.~{Alpay}, F.~{Colombo}, and I.~{Sabadini}.
\newblock {Schur functions and their realizations in the slice hyperholomorphic
  setting}.
\newblock To appear in {Integral Equations and Operator Theory}.

\bibitem{al_acap}
D.~Alpay and D.~Levanony.
\newblock Linear stochastic systems: a white noise approach.
\newblock {\em {Acta Applicandae Mathematicae}}, 110:545--572, 2010.

\bibitem{alp}
D.~Alpay, D.~Levanony, and A.~Pinhas.
\newblock Linear stochastic state space theory in the white noise space
  setting.
\newblock {\em {SIAM} {Journal of Control and Optimization}}, 48:5009--5027,
  2010.

\bibitem{MR2124899}
D.~Alpay, M.~Shapiro, and D.~Volok.
\newblock Rational hyperholomorphic functions in {$R^4$}.
\newblock {\em J. Funct. Anal.}, 221(1):122--149, 2005.


\bibitem{MR1129886}
Michael Artin.
\newblock {\em Algebra}.
\newblock Prentice Hall Inc., Englewood Cliffs, NJ, 1991.

\bibitem{bgr}
J.~Ball, I.~Gohberg, and L.~Rodman.
\newblock {\em Interpolation of rational matrix functions}, volume~45 of {\em
  Operator {T}heory: {A}dvances and {A}pplications}.
\newblock Birkh{\" a}user Verlag, Basel, 1990.

\bibitem{Bourbaki_Algebra}
N. Bourbaki.
\newblock {\em {Algebra. {I}. {C}hapters 1--3}}.
\newblock Springer-Verlag, 1989.

\bibitem{MR839186}
J.~W. Brewer, J.~W. Bunce, and F.~S. Van~Vleck.
\newblock {\em Linear systems over commutative rings}, volume 104 of {\em
  Lecture Notes in Pure and Applied Mathematics}.
\newblock Marcel Dekker Inc., New York, 1986.

\bibitem{MR1311923}
Paulo~D. Cordaro and Fran{\c{c}}ois Tr{\`e}ves.
\newblock {\em Hyperfunctions on hypo-analytic manifolds}, volume 136 of {\em
  Annals of Mathematics Studies}.
\newblock Princeton University Press, Princeton, NJ, 1994.

\bibitem{MR2584037}
Arthur~E. Frazho and Wisuwat Bhosri.
\newblock {\em An operator perspective on signals and systems}, volume 204 of
  {\em Operator Theory: Advances and Applications}.
\newblock Birkh\"auser Verlag, Basel, 2010.
\newblock Linear Operators and Linear Systems.

\bibitem{Fuhrmann}
P.A. Fuhrmann.
\newblock {\em Linear systems and operators in {H}ilbert space}.
\newblock McGraw-Hill international book company, 1981.

\bibitem{GS2_english}
I.M. Gelfand and G.E. Shilov.
\newblock {\em {Generalized functions. Volume 2}}.
\newblock Academic Press, 1968.

\bibitem{Groth55}
A.~Grothendieck.
\newblock {\em Produits tensoriels topologiques et espaces nucl{\'e}aires},
  volume~16.
\newblock Mem. Amer. Math. Soc, 1955.

\bibitem{MR1071708}
J.~{\'A}. Hermida~Alonso and T.~S{\'a}nchez-Giralda.
\newblock On the duality principle for linear dynamical systems over
  commutative rings.
\newblock {\em Linear Algebra Appl.}, 139:175--180, 1990.

\bibitem{Hida_Ikeda}
T.~Hida and N.~Ikeda.
\newblock Analysis on {H}ilbert space with reproducing kernel arising from
  multiple {W}iener integral.
\newblock In {\em Proc. {F}ifth {B}erkeley {S}ymp. {M}ath. {S}tat. {P}robab.
  {II}, part 1}, pages 117--143. {University of California Press}, 1967.

\bibitem{MR1502747}
Einar Hille.
\newblock A class of reciprocal functions.
\newblock {\em Ann. of Math. (2)}, 27(4):427--464, 1926.

\bibitem{Hille_Duke}
Einar Hille.
\newblock Contributions to the theory of {H}ermitian series.
\newblock {\em Duke Math. J.}, 5:875--936, 1939.

\bibitem{MR0000871}
Einar Hille.
\newblock Contributions to the theory of {H}ermitian series. {II}. {T}he
  representation problem.
\newblock {\em Trans. Amer. Math. Soc.}, 47:80--94, 1940.

\bibitem{MR1408433}
H.~Holden, B.~{\O}ksendal, J.~Ub{\o}e, and T.~Zhang.
\newblock {\em Stochastic partial differential equations}.
\newblock Probability and its Applications. Birkh\"auser Boston Inc., Boston,
  MA, 1996.

\bibitem{MR1851117}
Zhi-yuan Huang and Jia-an Yan.
\newblock {\em Introduction to infinite dimensional stochastic analysis},
  volume 502 of {\em Mathematics and its Applications}.
\newblock Kluwer Academic Publishers, Dordrecht, chinese edition, 2000.

\bibitem{MR0255260}
R.~E. Kalman, P.~L. Falb, and M.~A. Arbib.
\newblock {\em Topics in mathematical system theory}.
\newblock McGraw-Hill Book Co., New York, 1969.

\bibitem{MR751959}
Michael Reed and Barry Simon.
\newblock {\em Methods of modern mathematical physics. {I}}.
\newblock Academic Press Inc. [Harcourt Brace Jovanovich Publishers], New York,
  second edition, 1980.
\newblock Functional analysis.

\bibitem{MR1478165}
S.~Saitoh.
\newblock {\em Integral transforms, reproducing kernels and their
  applications}, volume 369 of {\em Pitman Research Notes in Mathematics
  Series}.
\newblock Longman, Harlow, 1997.

\bibitem{sansone}
G.~Sansone.
\newblock {\em Orthogonal functions}.
\newblock Dover {P}ublications, {I}nc., {N}ew--{Y}ork, 1991.
\newblock Revised {E}nglish {E}dition.

\bibitem{MR2105995}
Barry Simon.
\newblock {\em Functional integration and quantum physics}.
\newblock AMS Chelsea Publishing, Providence, RI, second edition, 2005.

\bibitem{SSR}
E.D. Sontag.
\newblock Linear systems over commutative rings: A survey.
\newblock {\em Ricerche di Automatica}, 7:1--34, 1976.

\bibitem{MR0372517}
G{\'a}bor Szeg{\H{o}}.
\newblock {\em Orthogonal polynomials}.
\newblock American Mathematical Society, Providence, R.I., fourth edition,
  1975.
\newblock American Mathematical Society, Colloquium Publications, Vol. XXIII.

\bibitem{Treves67}
F.~Treves.
\newblock {\em Topological vector spaces, distributions and kernels}.
\newblock Academic Press, 1967.

\bibitem{vage96b}
G.~V{\aa}ge.
\newblock A general existence and uniqueness theorem for {W}ick-{SDE}s in
  $(\mathcal S)^n_{-1,k}$.
\newblock {\em Stochastic {S}ochastic {R}ep.}, 58:259--284, 1996.

\bibitem{vage96}
G.~V{\aa}ge.
\newblock Hilbert space methods applied to stochastic partial differential
  equations.
\newblock In H.~K{\"o}rezlioglu, B.~{\O}ksendal, and A.S. {\"U}st{\"u}nel,
  editors, {\em Stochastic analysis and related topics}, pages 281--294.
  Birk{\"a}user, {B}oston, 1996.

\bibitem{zemanian}
A.H. Zemanian.
\newblock {\em Realizability theory for continuous linear systems}.
\newblock Dover {P}ublications, {I}nc., {N}ew--{Y}ork, 1995.

\bibitem{zhang}
T~Zhang.
\newblock Characterizations of white noise test functions and {H}ida
  distributions.
\newblock {\em Stochastics}, 41:71--87, 1992.
\end{thebibliography}
\def\cprime{$'$} \def\lfhook#1{\setbox0=\hbox{#1}{\ooalign{\hidewidth
  \lower1.5ex\hbox{'}\hidewidth\crcr\unhbox0}}} \def\cprime{$'$}
  \def\cprime{$'$} \def\cprime{$'$} \def\cprime{$'$} \def\cprime{$'$}

\end{document}